\input amstex
\documentstyle{amams} % input Annals of Mathematics macros.
\document
\annalsline{153}{2001}
\received{January 24, 2000}
\startingpage{767}
\font\tenrm=cmr10

\catcode`\@=11
\font\twelvemsb=msbm10 scaled 1100

%\font\ninemsb=msbm7 scaled 1100%msbm9
\font\ninemsb=msbm10 scaled 800
\newfam\msbfam
\textfont\msbfam=\twelvemsb  \scriptfont\msbfam=\ninemsb
  \scriptscriptfont\msbfam=\ninemsb
\def\msb@{\hexnumber@\msbfam}
\def\Bbb{\relax\ifmmode\let\next\Bbb@\else
 \def\next{\errmessage{Use \string\Bbb\space only in math
mode}}\fi\next}
\def\Bbb@#1{{\Bbb@@{#1}}}
\def\Bbb@@#1{\fam\msbfam#1}
\catcode`\@=12

 \catcode`\@=11
\font\twelveeuf=eufm10 scaled 1100
\font\teneuf=eufm10
\font\nineeuf=eufm7 scaled 1100%eufm9
\newfam\euffam
\textfont\euffam=\twelveeuf  \scriptfont\euffam=\teneuf
  \scriptscriptfont\euffam=\nineeuf
\def\euf@{\hexnumber@\euffam}
\def\frak{\relax\ifmmode\let\next\frak@\else
 \def\next{\errmessage{Use \string\frak\space only in math
mode}}\fi\next}
\def\frak@#1{{\frak@@{#1}}}
\def\frak@@#1{\fam\euffam#1}
\catcode`\@=12

%--------------- Author macros ---------------
\define \gln {\text{\rm GL}_n}
\define \gl {\text{\rm GL}}
\define \hofF {{\Cal H}(F)}

\define \br {\Bbb R}
\define \be {\Bbb E}
\define \bc {\Bbb C}
\define \bz {\Bbb Z}
\define \bq {\Bbb Q}
\define \ux {\bold x}
\define \uX {\bold X}
\define \uY {\bold Y}
\define \uo {\bold 0}
\define \um {\bold m}
\define \ue {\bold e}
\define \ua {\bold a}
\define \uj {\bold j}
\define \uy {\bold y}
\define \uz {\bold z}
\define \uL {\bold L}
\define \uK {\bold K}
%-------------- Author entries --------------------

\title{Decomposable form inequalities}
  \acknowledgements{Research partially supported by NSF grant DMS-9800859.}
\author{Jeffrey Lin Thunder}

  \institutions{Northern Illinois University, DeKalb, IL\\
{\eightpoint {\it E-mail address\/}: jthunder\@math.niu.edu}}

 \centerline{\bf Abstract}
\smallbreak
We consider Diophantine inequalities of the kind $|F(\ux )|\le m$,
where $F(\uX )\in \bz [\uX]$ is a homogeneous polynomial which
can be expressed as a product of $d$
homogeneous linear
forms in $n$ variables with complex coefficients and $m\ge 1$. 
We say such a form
is of finite type if the total volume of all real solutions to
this inequality is finite and if, for every $n'$-dimensional
subspace $S\subseteq \br ^n$ 
defined over $\bq$, the corresponding $n'$-dimensional volume
for $F$ restricted to $S$ is also finite.

We show that the 
number of integral solutions $\ux\in\bz ^n$ to our inequality above
is finite for all $m$ if and only if the form $F$ is of finite type.
When $F$ is of finite type, we show that the
number of integral solutions
is estimated asymptotically as $m\rightarrow \infty$
by the total volume of all real solutions. This generalizes a previous
result due to Mahler for the case $n=2$.
Further, we
prove a conjecture of W. M.\ Schmidt, showing that  for $F$ of finite type
the number of integral solutions is bounded above by $c(n,d)m^{n/d},$
where $c(n,d)$ is an effectively computable constant depending only
on $n$ and $d$.

\bigbreak
\intro

In this paper we consider forms in $n>1$ variables of the type
$F(\uX)=\prod _{i=1}^dL_i(\uX)\in\bz [\uX ],$ 
where each $L_i(\uX)\in\bc [\uX ]$ is
a linear form.
For a positive integer $m$
we are interested in the integer
solutions $\ux\in\bz ^n$ to the inequality
$$|F(\ux )|\le m.\tag 1$$

Consider the case when $n=2$, $d>n$ and $F(\uX)$ is
irreducible over $\bq$. 
Thue's famous result in [T] is that the number of integer solutions to
(1) in this case is finite.
Later, Mahler in [M] estimated the
number $N_F(m)$ of\break such solutions as follows. Let $A(F)$ denote
the area of the planar region\break
$\{\ux\in\br ^2\: |F(\ux )|\le 1\},$
so that $m^{2/d}A(F)$ is the measure of the set of $\ux\in\br ^2$
that satisfy (1). (The hypothesis that $F$ is irreducible forces the
discriminant to be nonzero, which
implies that this area is finite.) Then 
$$\left |N_F(m)-m^{2/d}A(F)\right |=O\big (m^{1/(d-1)}\big )$$
as $m\rightarrow \infty$,
where the implicit constants depend on $d$ and $F$.
We also have a result due
to Schmidt [S3, Chap.~III, Theorem 1C] which states that for irreducible $F$,
$N_F(m)\ll dm^{2/d}(1+\log m^{1/d})$
with an absolute implicit constant.

Other than results for the case $n=2$, little has been
published on this question. Ramachandra in [R] proved that for
norm forms of the type $N_{K/\bq}(X_1+\alpha X_2+\cdots +\alpha ^{n-1}
X_n)$, where $K=\bq (\alpha )$ is a number field of degree $d\ge 8n^6$ and
$N_{K/\bq}$ denotes the norm from $K$ to $\bq$,
one has
$$|N_F(m)-m^{n/d}V(F)|=O(m^{\varepsilon +(n-1)/(d-n+2)})$$
for any $\varepsilon >0$ as $m\rightarrow\infty,$ 
where the implicit constant depends on both
$F$ and $\varepsilon$, and $V(F)$ denotes the volume analogous to the
area $A(F)$ above. Note that by the homogeneity of $F$, $m^{n/d}V(F)$ is
the volume of the set of all real solutions to (1).
Of course,
one needs the subspace theorem to approach the general case. For norm
forms,
Schmidt showed in [S1] that the number of solutions to (1) is finite for
all $m$ if
and only if $F$ is a nondegenerate.
Evertse has shown in [E3] that for nondegenerate norm forms $F$ of
degree $d$ in $n$ variables, one has
$$N_F(m)\le (16d)^{{1\over 3}(n+7)^3}m^{(n+\sum _{i=2}^{n-1}i^{-1})/d}\times
(1+\log m)^{{1\over 2}n(n+1)}.$$

The results above are of two different flavors. On the one hand
the natural heuristic is that, in the absence of a compelling reason
to the contrary, one expects that the volume of the region in $\br ^n$ 
defined by (1) should approximate the number of integral solutions to (1).
The results of Mahler and Ramachandra above verify this in special cases.
On the other hand, when $N_F(m)$ is finite one expects that it should be
bounded above by a function independent of the specific
coefficients of $F$. This was proven by Evertse in [E1] for the case
$n=2$, and another result of Schmidt in [S2] confirms this 
in the general case of products of nondegenerate norm forms. 
Schmidt's absolute upper bound above in the case $n=2$ 
appears to be the right order of 
magnitude in terms of $m$ (up to the logarithmic term).
In fact, in [S2] Schmidt makes the conjecture that
$N_F(m)\ll m^{n/d}$ for all nondegenerate norm forms of degree $d$ in
$n$ variables, where the implicit
constant depends only on $n$ and $d$.  Evertse's result above comes
close to this. 

When one tries to reconcile the heuristic with Schmidt's conjecture, one
is led to the conjecture that $V(F)\ll 1$ 
for nondegenerate norm forms. 
This was shown to be true in [B] for
the case $n=2$, and was shown to be true for 
general forms in $d>n$ variables with nonzero discriminant in [BT]. 
 
Returning to our heuristic, what would be a ``compelling reason"
for $N_F(m)$ to {\it not} be approximated by the volume? One such
reason comes immediately to mind. It is typically the case that, though
the volume $V(F)$ may be finite, the lower
dimensional volume of the region defined by (1) cut by a 
hyperplane is infinite. 
If such a hyperplane were defined over $\bq,$
then that rational hyperplane might contain
infinitely many integral points. With this in mind, we say $F$ is
of {\it finite type} if $V(F)$ is finite, and the same is true for
$F$ restricted to any nontrivial rational subspace. Note in particular
that if $F$ is of finite type, it does not vanish at any nonzero
rational point. When $F$ is of finite type, then, we rule out this
``compelling reason." Since $N_F(m)$ can be infinite if $F$ is a degenerate
norm form, this could be a ``compelling reason" as well. But degeneracy
of a norm form is a rather algebraic concept, and it is not immediately
clear what the connection is between this and
the more geometric concept of the volume $V(F)$.
 
The purpose of this paper is to answer the following questions: When is
$V(F)$ finite? More correctly, can one determine rather simply from a given
factorization of $F$ whether $V(F)$ is finite or not? If $V(F)$ is finite, is
$V(F)\ll 1$?
When is $N_F(m)$ finite for all $m$? If $N_F(m)$ is finite, is it approximated
by $m^{n/d}V(F)$? If $N_F(m)$ is finite, is $N_F(m)\ll m^{n/d}$?
We will  prove Schmidt's conjecture
and more.
Here and from now on, all implicit constants in the
$\ll$ notation depend only (and explicitly) on $n$ and $d$.

\nonumproclaim{Theorem 1} Let $F$ be a decomposable
form of degree $d$ in $n$ variables
with integral coefficients{\rm .} If $V(F)$ is finite and $F$ does not vanish
at a nonzero integral point{\rm ,}
then
$V(F)\ll 1${\rm . }
\endproclaim

\nonumproclaim{Theorem 2} Let $F$ be a decomposable form of degree $d$ in
$n$ variables with integral coefficients{\rm .} Then $N_F(m)$ is finite
for all $m$
if and only if $F$ is of finite type{\rm .} If $F$ is of finite type{\rm ,} then 
$N_F(m)\ll m^{n/d}.$
\endproclaim

Apparently nondegenerate
norm forms are of finite type. This could be shown more directly, though
it is not a simple consequence of the definition of nondegenerate. 
The answer to our question regarding the finiteness of $V(F)$ 
requires further notation, so we leave it for
the next section (see the proposition below). We only
remark here that it
is necessary that $d>n$ in order for $V(F)$ to be finite except for
the case of a positive definite quadratic form in two variables.

\nonumproclaim{Theorem 3} Let $F$ be decomposable form of degree $d$ in
$n$ variables  with integral
coefficients{\rm .} If $F$ is of finite type{\rm ,} then
there are $a(F),\ c(F)\in\bq$ satisfying
$$1\le a(F)\le{d\over n}-{1\over n(n-1)}$$
and 
$${(d-n)\over d}\le c(F)< {d\choose n}(d-n+1)$$
such that
$$|N_F(m)-m^{n/d}V(F)|\ll m^{(n-1)/(d-a(F))}(1+\log m)^{n-2}
\hofF ^{c(F)}.$$
If the
discriminant is not zero{\rm ,} then we may take $a(F)=1$ and 
$c(F)=\break {d-1\choose n-1}-1.$
\endproclaim

The quantities $a(F),\ c(F)$ and $\hofF$ appearing in Theorem 3
are explicitly defined in the next section. 
Note that $(n-1)/(d-a(F))<n/d$ in Theorem 3, so that the estimate
for $N_F(m)$ given is not trivial. Theorem 3 is a broad generalization of
Mahler's result above. 

This paper is organized as follows. Section 1 introduces some notation
and defines some quantities connected to $F$ which will be used 
throughout. In the next section we derive some general results
concerning the height $\hofF$. Sections 3 and 4 are the technical heart
of the paper where we see that solutions to (1) lie in subsets
of certain convex regions (these regions are parallelopipeds 
if $F$ factors over $\br$) and we garner pertinent information
about these convex regions. Section 5 deals with the case when $V(F)$ is
infinite. The next two sections are devoted
to estimating volumes connected with (1) and analyzing the set of
integral solutions 
to (1).
The proofs of our theorems follow in the last section, using an inductive
argument on the number of variables $n$.

\section{Definitions and a linear programming result}

Throughout the rest of this paper,
$F(\uX )=\prod _{i=1}^dL_i(\uX )\in\bz [\uX ]$
will denote a decomposable form of degree $d$ in $n$ variables with integral
coefficients
and $m\ge 1$ will be a fixed real number. The case where $F$ is 
a power of a positive definite quadratic form in two variables is exceptional
and our questions posed in the introduction are trivially answered in this
case, so from now on we will assume that $F$ is not such a form.

We will use the notion of ``equivalent forms."
If $F$ is a decomposable form in $n$ variables and $T\in\gln (\bz)$, then we can
compose $F$ with $T$ to get a new form $G(\uX )=F\circ T(\uX )$.
Since $\det (T)=\pm 1,$ we have $V(F)=V(G)$. 
Further, the integral solutions to (1)
are in one-to-one correspondence (via $T^{-1}$) with the integral
solutions to $|G(\ux)|\le m$. Because of this, we say two forms
$F$ and $G$ are equivalent if there is a $T\in\gln (\bz )$ with
$G=F\circ T$. The freedom to choose a representative from each
equivalence class will be used to our advantage.

We now proceed with some definitions and notation.
We will denote the usual $L_2$ norm of $\ux\in\bc ^n$
by $\| \ux \|$. We will denote the coefficient vector of a linear
form $L_i(\uX )$ by
$\uL _i\in\bc ^n$.  Complex conjugation will be denoted by
an overline: $\overline{\alpha}$. This notation will be extended
to vectors as well, e.g., $\overline{\uL}.$ Elements of $\bc ^n$ will be
viewed as $1\times n$ matrices (i.e., row vectors) and a
superscript $^{tr}$ will denote the transpose of a matrix, so
that $\uL ^{tr}$ is a column vector for a coefficient vector $\uL$.

We define the {\it height} of $F$ to be
\smallbreak
\centerline{${\displaystyle \hofF :=\prod _{i=1}^d\|\uL _i\|.}$}
\smallbreak\noindent 
Note that $\hofF$ is
actually independent of the particular factorization of $F$ used,
though it is not preserved under equivalence.

Given a factorization of $F$, let $I(F)$ denote the set of
all ordered $n$-tuples $(\uL _{i_1},\ldots ,\uL _{i_n})$ of linearly
independent coefficient vectors.
We let $b(\uL _i)$ denote the number of $n$-tuples in $I(F)$ where
$\uL _i$ occurs and let $b(F)$ denote the maximum of these $b(\uL _i)$.
Note that $b(F)$ is 
preserved under equivalence and is independent of the factorization used. 
Let $I'(F)\subset I(F)$ denote those $n$-tuples with $i_1<i_2<\cdots <i_n.$
Letting $|\cdot |$ denote the cardinality, we have
$$|I(F)|=n!|I'(F)|\le n!{d\choose n},\tag 2$$
with equality if and only if the discriminant of $F$ is not zero.

Let $J(F)$ be the subset of $I(F)$ consisting of $n$-tuples 
that satisfy
the following restriction: if $j<n,$ then either $\uL _{i_{j+1}}$ is 
proportional to $\overline{
\uL _{i_j}}$ or $\overline{\uL _{i_j}}$ is in the span of 
$\uL _{i_1},\ldots ,\uL _{i_j}.$ 
If $J(F)$ is not empty, we let
$$a(F)=\max \left \{ \text{the number of $\uL _i$ in the span of
$\uL _{i_1},\ldots ,\uL _{i_j}$}\over j\right \},$$
where the maximum is over all $n$-tuples in $J(F)$ and $j=1,\ldots ,n-1.$
If $J(F)$ is empty, we leave $a(F)$ undefined.
Note that the number of factors in the span of $\uL _{i_1},\ldots ,\uL _{i_n}$
is $d$ for all $n$-tuples in $I(F)$.
We will see later
(see Lemma 5 below) that $J(F)$ is in fact empty only when $I(F)$ is.
Note that $a(F)\ge 1$ if it is defined, with equality if and only
if the discriminant of $F$ is not zero.

We can now state our characterization of finite volume in terms of
the factorization of $F$.

\nonumproclaim{Proposition} For a decomposable form $F$ as above{\rm ,}
$V(F)$ is finite if and only if $a(F)$ is defined and less than $d/n${\rm .}
\endproclaim

The proposition will be proven in Section 7 below. 

We now continue with
some definitions.
Let
$$c(F)=\cases {d-1\choose n-1}-1&\matrix \text{if the discriminant}\hfill\\ \noalign{\vskip-6pt}
\text{of $F$ is not zero,}\hfill\endmatrix
\\
{b(F)\over n!a(F)}\big (d-(n-1)a(F)\big ) -{1\over a(F)}&
\text{otherwise,}\endcases$$
whenever $a(F)$ is defined.
This quantity occurs as an exponent on $\hofF$ in our arguments; we give
it a name for notational convenience.

The {\it semi-discriminant} of $F$, which we denote by $S(F)$, is given by
$$S(F):= \prod 
\det (\uL ^{tr}_{i_1},\ldots ,\uL _{i_n}^{tr}),$$
where the product is over all $n$-tuples in $I(F)$ when $I(F)$ is not empty, 
and $S(F)=0$ otherwise. 
Unlike $\hofF$, the semi-discriminant can be dependent on the factorization.
If 
$$F(\uX )=\prod _{i=1}^dL _i(\uX )=\prod _{i=1}^d\alpha _iL_i(\uX )$$ are
two different factorizations of $F$, then the semi-discriminant for the first
will equal that for the second if and only if
$$\prod _{i=1}^d\alpha _i^{b(\uL _i)}=1.$$
Hence, the semi-discriminant is independent of the factorization if and only
if $b(\uL _i)=b(F)$ for all $i$.
This is not always the case, as the example $F(\uX )=X_1^2X_2X_3\cdots X_n$
shows. 
To deal with this nonuniqueness, we introduce a quantity which we
call the {\it normalized semi-discriminant}, denoted by $NS(F)$ and defined
by
$$NS(F):=\prod {
\det (\uL ^{tr}_{i_1},\ldots ,\uL _{i_n}^{tr})\over
\|\uL _{i_1}\|\cdots\|\uL _{i_n}\|}= {S(F)\over \|\uL _1\|^{b(\uL _1)}\cdots
\|\uL _d\|^{b(\uL _i)}},$$
where the product is over all $n$-tuples in $I(F)$. Then $|NS(F)|$ is entirely
determined by the form $F$.
It is not preserved under
equivalence.

We end this section with a simple linear programming result which will
be needed later.

\nonumproclaim{Lemma 1} Let $k$ be a positive integer{\rm .} Let $b_1\le \cdots \le b_k$
be a nondecreasing sequence of real numbers and let $A>0${\rm .} Then the minimum
value of\break $x_1b_1+\cdots +x_kb_k$ subject to the restrictions
$$\aligned  x_i&\ge 0\hskip8pt \qquad\text{all $i$},\\
x_1+\cdots +x_j&\le jA \qquad\text{all $j$},\\
x_1+\cdots +x_k&=kA,\endaligned $$
is achieved when $x_i=A$ for all $i$.\endproclaim 

\demo{{P}roof} We prove this by induction on $k$. The case $k=1$
is trivial, so assume $k>1$. 

Suppose $x_1,\ldots ,x_k$ satisfy the restrictions given. Let $i$ be minimal
such that $x_i>0$. If $i>1,$ then 
$$x_j'=\cases x_j&\text{if $j\neq i, i-1$},\\ x_i/2&\text{otherwise}\endcases$$
also satisfy the restrictions, and $x_1b_1+\cdots +x_nb_n\ge x_1'b_1+\cdots
+x_n'b_n$ since $b_{i-1}\le b_i$. This shows that there is a solution to 
our 
problem where $x_1>0$. On the other hand, it is well known 
that any solution to such a problem
occurs at a vertex of the convex region determined by the restrictions. Such
a vertex has $x_1=0$ or $A$, so the minimum can be achieved when $x_1=A.$

We now invoke the induction hypothesis, which says that the minimum value
of $x_2b_2+\cdots +x_kb_k$ subject to the restrictions
$$\aligned  x_i&\ge 0 \phantom{j-1)A} \qquad\text{all $i>1$},\\
x_2+\cdots +x_j&\le (j-1)A\, \qquad\text{all $j>1$},\\
x_2+\cdots +x_k&=(k-1)A,\endaligned $$
is achieved when $x_i=A$ for all $i>1$.
\enddemo  

\section{Inequalities involving the height}

\nonumproclaim{Lemma 2} For any factor $L_i(\uX )$ of $F(\uX )${\rm ,} $\uL _i$ is
proportional to a vector $\uL _i'$ with algebraic coefficients in a number field
of degree no greater than $d${\rm ,} and the field height $H(\uL _i')$ satisfies
$H(\uL _i')\le \hofF.$ In particular{\rm ,} $\hofF\ge 1${\rm .}
\endproclaim

See [S3] for a definition of $H(\uL )$. This is the usual field height
(not absolute height) using $L_2$ norms at the infinite places.

\demo{{P}roof} 
Suppose first that $F$ is irreducible over $\bq$. 
It is known that $F(\uX )=aN_{K/\bq}\big (L(\uX )\big )$, where $a$ is a
nonzero rational number, $K$ is a number field of degree equal to the degree
of $F$ and $N_{K/\bq}$ denotes the norm from $K$ to $\bq$. Thus, any factor
of $F$ is proportional to some conjugate of $L(\uX )$. The coefficient
vectors of these conjugates
all have the same field height (see the remark on p.~23 of [S3]).
Further, by [S3 Chap.~III, Lemma 2A],
$\hofF={\rm cont}(F)H(\uL )$, where ${\rm cont}(F)$ denotes the content of $F$. 
Since the content of $F$ is a positive integer, we get $H(\uL )\le \hofF$.
Since the field height function $H\ge 1$, the lemma is true when 
$F$ is irreducible over $\bq$.

In general,
$$F(\uX )=\prod _{l=1}^kF_l(\uX ),$$
where each $F_l$ is a form with integral coefficients which is
irreducible over $\bq$. Any linear factor $L_i(\uX )$ of $F$ is a factor of
some $F_{l_i}$. By what we have shown, $\uL _i$ is proportional to
an $\uL _i'$ with algebraic coefficients in a number field of degree
no greater than the degree of $F_{l_i}$ and satisfying $H(\uL _i')\le {\Cal H}
(F_{l_i})$. The degree of $F_{l_i}$ is certainly no larger than the 
degree of $F$, and
$$\hofF =\prod _{l=1}^k{\Cal H}(F_l).$$ 
We have shown that ${\Cal H}(F_l)\ge 1$ for all $l$, so 
$\hofF\ge{\Cal H}(F_{l_i})$ and the lemma is proven.
\enddemo

\nonumproclaim{Lemma 3} If $I(F)$ is not empty{\rm ,} then
$$|NS(F)|\ge\hofF ^{-b(F)}.$$ 
For any $n$-tuple in $I(F)$ we have
$${|\det (\uL _{i_1}^{tr},\ldots ,\uL _{i_n}^{tr})|\over\prod _{j=1}^n
\|\uL _{i_j}\|}\ge\hofF ^{-b(F)/n!}.\tag 3$$
\endproclaim

\demo{{P}roof} Since $|NS(F)|$ is independent of
the factorization used, we may choose any one we wish. First factor $F$ into
a product of forms with integral coefficients which are irreducible over
$\bq$,
$$F(\uX )=\prod _{l=1}^kF_l(\uX ),$$
as in the proof of Lemma 2 above.  Write each $F_l(\uX )$ as a rational 
multiple of a norm form as above in the proof of Lemma 2.

Since $F$ has rational coefficients, it is invariant under any
element $\sigma$ of the Galois group of $\overline{\bq}$ over $\bq$, where
$\overline{\bq}\subset\bc$ is the algebraic closure of $\bq$ in $\bc$.
Thus, any element of the Galois group must
take our factorization of $F$ to another, say
$$\sigma (\uL _i)=\beta _i\uL _{\sigma '(i)},$$
where $\beta _i\in\bc ^{\times}$ and $\sigma '$ is
an element of the permutation group of $\{1,\ldots ,d\}$. 
Also, $\sigma (S(F))$ is equal to the semi-discriminant with this
factorization given by $\sigma$. 
For any $n$-tuple  $(\uL _{i_1},\ldots ,\uL _{i_n})\in I(F)$, 
$$0\neq \sigma \big (\det (\uL _{i_1}^{tr},\ldots ,\uL _{i_n}^{tr})\big )=
\det (\uL _{\sigma '(i_1)}^{tr},\ldots ,\uL _{\sigma '(i_n)}^{tr})\times
\prod _{j=1}^n\beta _{i_j}.$$
Since $\sigma '$ is a permutation, in this manner we
see that $b(\uL _i)=b(\uL _{\sigma '(i)}
)$ for any $i$. But the Galois group acts transitively on the factors of
norm forms, so we conclude that 
$b(\uL _i)=b(\uL _j)$
whenever $L_i(\uX )$ and $L_j(\uX )$ are factors of the same irreducible 
$F_l(\uX )$, i.e.,
all the linear factors of a given $F_l$ have the same $b$ value.
Let $b_l$ denote the $b$ value of the linear factors of $F_l$ for each
$l=1,\ldots ,k$.

Suppose $\uL _{i_1},\ldots ,\uL _{i_{d'}}$ are the coefficient vectors
of the linear factors of some
$F_l$. Just like $F(\uX )$, $F_l(\uX )$ has integral coefficients and 
is invariant under $\sigma$. Hence,
$$\prod_{j=1}^{d'} \beta _{i_j}=1
=\prod _{j=1}^{d'}\beta _{i_j}^{b_l}.$$
Taking into account the different factors $F_l$ of $F$, we are led to
$$\prod _{i=1}^d\beta _i^{b(\uL _{\sigma '(i)})}=1.$$
As remarked in Section 1, this shows that
the semi-discriminant $S(F)$ is the same for our initial factorization of $F$
and the factorization induced by $\sigma$. So our $S(F)$ is invariant
under the Galois group of $\overline{\bq}$, and hence a rational number.
It is nonzero since $I(F)$ is not empty.

Let $v$ be any place of $\overline{\bq}$. Then Hadamard's
inequality gives
$${|\det (\uL _{i_1}^{tr},\ldots ,\uL _{i_n}^{tr})|_v\over\|\uL _{i_1}\|_v
\cdots \|\uL _{i_n}\|_v}\le 1,$$
where $\|\cdot\| _v$ denotes the usual $L_2$ norm if $v|\infty$ and
the sup norm otherwise. 
In particular,
$$|S(F)|_v\le \prod _{i=1}^d\|\uL _i\| ^{b(\uL _i)} _v.\tag 4$$
By the definition of $b_l$ we have
$$\prod _{i=1}^dL _i(\uX )^{b(\uL _i)}=F_1(\uX )^{b_1}\cdots F_k(\uX )^{b_k}.
\tag 5$$
We let ${\bold F_l}$  denote the coefficient vector of $F_l$ for each
$l=1,\ldots ,k$. If $v$
is non-archimedean, then Gauss' 
lemma together with (4) and (5) gives
$$|S(F)|_v\le\prod _{i=1}^d\|\uL _i\| _v^{b(\uL _i)}=\|{\bold F}_1\| _v^{b_1}
\cdots \|{\bold F}_k\|_v ^{b_k}\le 1.$$
This holds for any non-archimedean place, so $|S(F)|$ is a positive integer.
In particular, $|S(F)|\ge 1$.
By Lemma 3, ${\Cal H}(F_l)\ge 1$ for all $l$, so that by (5)
$$\prod _{i=1}^d\|\uL _i\|^{b(\uL _i)}={\Cal H}(F_1)^{b_1}\cdots {\Cal H}(F_k
)^{b_k}\le
{\Cal H}(F_1)^{b(F)}\cdots {\Cal H}(F_l)^{b(F)}=\hofF ^{b(F)}.$$
Hence $|NS(F)|\ge\hofF ^{-b(F)}$ with this factorization of $F$. 

As for (3), we note that
$$|NS(F)|^{1/n!}= \prod
{|\det (\uL _{i_1}^{tr},\ldots ,\uL _{i_n}^{tr})|\over\prod _{j=1}^n
\|\uL _{i_j}\|}\ge\hofF ^{-b(F)/n!},$$
where the product is over all $n$-tuples of $I'(F)$. We saw above that
each factor in this middle product is no greater than 1, thus each
factor is 
bounded below by our lower bound for
$|NS(F)|^{1/n!}.$
\enddemo

\section{Bounds for linear factors}

In this section our goal  is to show that for any solution $\ux\in\br ^n$ of (1),
there is an $n$-tuple in $I(F)$ with the product
$|L_{i_1}(\ux )|\cdots |L_{i_n}(\ux )|$ relatively small. We start with
a general result which says that $n$ linearly independent linear forms
cannot simultaneously be small at $\ux$.

\nonumproclaim{Lemma 4} Let $\ux\in\br ^n\setminus \{\uo \}$ and let $
L_1(\uX ),\ldots ,L_n(\uX )$ be $n$ linearly independent linear forms{\rm .}
Suppose that
$${|L_j(\ux )|\over \|\uL _j\|}\ge {|L_i(\ux )|\over \|\uL _i\|}$$
for $i=1,\ldots ,n${\rm .} Then
$${|L _j(\ux )|\over\|\uL _j\|}
\ge {\|\ux\| |\det (\uL ^{tr}_1,\ldots ,\uL ^{tr}_n)|
\over n^{n/2}\prod _{i=1}^n\|\uL _i\|}.$$
\endproclaim

\demo{{P}roof} Without loss of generality we may assume $\| \uL _i\| =1$
for all $i$ and $\|\ux\| =1$. Let $T$ denote the $n\times n$
matrix with rows $\uL _i$ and write
$$\aligned {\frak m}&=\min _{\| \uy \| =1}\left \{\| T \uy ^{tr} \| \right \}
\qquad\text{and}\\
{\frak M}&=\max _{\|\uy\| =1}\left \{\|T \uy ^{tr}\|\right \}.\endaligned$$

Suppose $\| T\ux _1^{tr}\|={\frak m}$ and $\|\ux _1\| =1.$ 
Choose $\ux _2,\ldots ,
\ux _n\in \br ^n$, all of length~1, that also satisfy
 $|\det (\ux _1^{tr}, \ldots
,\ux _n^{tr})|=1.$ We then have
$$\aligned |\det (T)|=|\det (T)||\det (\ux _1^{tr},\ldots ,\ux _n^{tr})|&=
|\det (T\ux _1^{tr},\ldots ,T\ux _n^{tr})|\\
&\le\prod _{l=1}^n\|T\ux _l^{tr}\|\\
&\le {\frak m}{\frak M}^{n-1}.\endaligned$$

Since $\| \uL _i\| =1$ for all $i$ we have ${\frak M}\le \sqrt{n}$, so that
$${\frak m}\ge n^{(1-n)/2}|\det (T)|.$$
By the hypothesis, $|L_j(\ux )|\ge |L_i(\ux )|$ for all $i$,
so that
$$\sqrt{n}|L_j(\ux )|\ge\|T\ux ^{tr}\|\ge {\frak m}.$$
Combining these last two inequalities yields the lemma.
\enddemo

\nonumproclaim{Lemma 5} 
Suppose $I(F)$ is not empty{\rm .} Then
$a(F)$ is defined{\rm .} If $a(F)<d/n${\rm ,} then 
for every $\ux\in\br^n$
there is an $n$\/{\rm -}\/tuple in $J(F)$ such that
$${\prod _{j=1}^n|L_{i_j}(\ux )|\over|\det (\uL _{i_1}^{tr},\ldots ,
\uL _{i_n}^{tr})|}\ll \left ({|F(\ux )|\over
\|\ux \|^{d-na(F)}}\right )^{1/a(F)}\hofF^{c(F)}.\tag 6$$
\endproclaim

\demo{{P}roof} Suppose $I(F)$ is not empty and let $\ux\in\br ^n.$
We define minima $\lambda _1\le\lambda _2\le\cdots\le\lambda
_n$ and choose indices $i_1,\ldots ,i_n$ as follows.
Let
$$\lambda _1=\min\{|L_i(\ux )|/\|\uL _i\|\},$$
where the minimum is over all factors $L_i(\uX )$ of $F(\uX )$. 
Choose $i_1$
such that 
$$|L_{i_1}(\ux )|/\|\uL _{i_1}\|=\lambda _1.$$
 We then continue
recursively, letting
$$\lambda _{j+1}=\min\{|L_i(\ux )|/\|\uL _i\|\}\ge \lambda _j,$$
where the minimum is over all factors $L_i(\uX )$ where $\uL _i$ is
not in the span of $\uL_{i_1},\ldots ,
\uL_{i_j}$, for $j=1,\ldots ,n-1$. We choose $i_{j+1}$ such that 
$\uL _{i_{j+1}}$
is not in the span of $\uL _{i_1},\ldots ,\uL _{i_j}$ and
$$|L_{i_{j+1}}(\ux )|/\|\uL _{i_{j+1}}\|=\lambda _{j+1},$$
 with the stipulation
that $\uL_{i_{j+1}}$ is proportional to $\overline{\uL _{i_j}}$ if 
$\overline{\uL _{i_j}}$ is
not in the span of $\uL _{i_1},\ldots ,\uL _{i_j}$. (Note that if this
were the case, then $\lambda _{j+1}=\lambda _j$, so that such a choice for
$i_{j+1}$ is possible.) These minima are well defined since $I(F)$ is not
empty, implying that the  set of all $\uL _i$ has rank $n$. By
construction, $(\uL _{i_1},\ldots ,\uL _{i_n})\in J(F)$, so $a(F)$ is defined.

Now suppose $a(F)<d/n$ and $\ux\in\br ^n$. If $F(\ux )=0,$
then (6) trivially holds since $\lambda _1=0$. So we may as well assume
$F(\ux )\neq 0$, which implies that $\lambda _1>0$. Let $a_1$ be the
number of $\uL _i$ which are linearly dependent on $\uL _{i_1}$.
For $j>1$ let $a_j$ be the number of $\uL _i$ which are in the span of
$\uL _{i_1},\ldots ,\uL _{i_j}$ but not in the span of $\uL _{i_1},\ldots ,
\uL _{i_{j-1}}$. If $\uL _i$ is in the span of 
$\uL _{i_1},\ldots ,\uL _{i_j}$ but not in the span of $\uL _{i_1},\ldots ,
\uL _{i_{j-1}}$, then $|L_i(\ux )|/\|\uL _i\|\ge\lambda _j$ by the
definition of $\lambda _j$. Thus,
$${|F(\ux )|\over \hofF}=\prod _{i=1}^d{|L_i(\ux )|\over \|\uL _i\|}\ge
\prod _{j=1}^n\lambda _j^{a_j}.\tag 7$$

By definition, $a_1+\cdots +a_j$ is the number of $\uL _i$ in the
span of $\uL _{i_1},\ldots ,\uL _{i_j}$ for $1\le j\le n$. This implies that
$$a_1+\cdots +a_j\le ja(F)\qquad 1\le j<n\tag 8$$ 
and $a_1+\cdots +a_n=d$. 
Let $s=a_1+\cdots +a_{n-1}$, so $s\le (n-1)a(F)$ by (8).
Since $\lambda _n\ge \lambda _{n-1},$ we see that
$$\lambda _n^{a_n}\lambda _{n-1}^{a_{n-1}}\ge\lambda _n^{a_n-((n-1)a(F)-s)}
\lambda _{n-1}^{a_{n-1}+((n-1)a(F)-s)}.$$
Define $a_j'$ by
$$a_j'=\cases a_n-((n-1)a(F)-s)=d-(n-1)a(F)&\text{for $j=n$,}\\
a_{n-1}+((n-1)a(F)-s)&\text{for $j=n-1$,}\\
a_j&\text{otherwise.}\endcases$$
Then (7) and (8) hold with $a_j'$ in place of $a_j$, and also
$$a_1'+\cdots +a_{n-1}'=(n-1)a(F).\tag 9$$

Because of (8) and (9), we can use Lemma 1 with $k=n-1$,
$b_j=\log\lambda _j$ and $A=a(F)$. We get
$$\prod _{j=1}^{n-1}\lambda _j^{a_j'}\ge\prod _{j=1}^{n-1}\lambda _j^{a(F)}.$$
This and (7) imply that
$$\aligned {|F(\ux )|\over \hofF}=\prod _{i=1}^d{|L_i(\ux )|\over \|\uL _i\|}&
\ge
\prod _{j=1}^n\lambda _j^{a_j'}\\
&=
\lambda _n^{d-na(F)+a(F)}\prod _{j=1}^{n-1}\lambda _j^{a_j'}\\
&\ge
\lambda _n^{d-na(F)}\prod _{j=1}^n\lambda _j^{a(F)}\\
&=\lambda _n^{d-na(F)}
\left (\prod _{j=1}^n{|L_{i_j}(\ux )|\over \|\uL _{i_j}\|}\right )^{a(F)}.
\endaligned$$
By Lemma 4,
$$\lambda _n={|L_{i_n}(\ux )|\over \|\uL _{i_n}\|}\gg \|\ux\|
{|\det (\uL _{i_1}^{tr},\ldots ,\uL _{i_n}^{tr})|\over
\prod _{j=1}^n\|\uL _{i_j}\|}.$$
So by (3)
$$\aligned {|F(\ux )|\over \hofF}&\ge
\lambda _n^{d-na(F)}
\left (\prod _{j=1}^n{|L_{i_j}(\ux )|\over \|\uL _{i_j}\|}\right )^{a(F)}\\
&\gg \|\ux\|^{d-na(F)}
\left ({|\det (\uL _{i_1}^{tr},\ldots ,\uL _{i_n}^{tr})|\over
\prod _{j=1}^n\|\uL _{i_j}\|}\right )^{d-na(F)}
\left (\prod _{j=1}^n{|L_{i_j}(\ux )|\over \|\uL _{i_j}\|}\right )^{a(F)}\\
&= \|\ux\|^{d-na(F)}
\left ({|\det (\uL _{i_1}^{tr},\ldots ,\uL _{i_n}^{tr})|\over
\prod _{j=1}^n\|\uL _{i_j}\|}\right )^{d-(n-1)a(F)} \\
&\qquad \times\ 
\left ({\prod _{j=1}^n|L_{i_j}(\ux )|\over 
|\det (\uL _{i_1}^{tr},\ldots ,\uL _{i_n}^{tr})|}
\right )^{a(F)}\\
&\ge \|\ux \|^{d-na(F)}\hofF ^{-b(F)\big (d-(n-1)a(F)\big )/n!}
\left ({\prod _{j=1}^n|L_{i_j}(\ux )|\over 
|\det (\uL _{i_1}^{tr},\ldots ,\uL _{i_n}^{tr})|}
\right )^{a(F)}.\endaligned$$
This proves (6) in the case where the discriminant is zero.

When the discriminant is not zero
we can do somewhat better. First of all, we have $a(F)=1$. Letting
$L_{i_1},\ldots ,L _{i_n}$ be as above, we see that
$${|F(\ux )|\over \hofF}\ge \prod _{l\neq i_1,\ldots ,i_n}
 {|L_l(\ux )|\over \|\uL _l\|}\times
\prod _{j=1}^n{|L_{i_j}(\ux )|\over \|\uL _{i_j}\|}.$$
By Lemma 4, 
$$\aligned {|F(\ux )|\over \hofF}&\gg
\prod _{l\neq i_1,\ldots ,i_n}
{\|\ux\|\cdot |\det (\uL _l^{tr},\uL _{i_1}^{tr},\ldots ,\uL _{i_{n-1}}^{tr})|
\over\|\uL _l\|\cdot \|\uL _{i_1}\|\cdots \|\uL _{i_{n-1}}\|}
\prod _{j=1}^n{|L_{i_j}(\ux )|\over \|\uL _{i_j}\|}\\
&=\|\ux\| ^{d-n}
\prod _{l\neq i_1,\ldots ,i_{n-1}}
{|\det (\uL _l^{tr},\uL _{i_1}^{tr},\ldots ,\uL _{i_{n-1}}^{tr})|
\over\|\uL _l\|\cdot \|\uL _{i_1}\|\cdots \|\uL _{i_{n-1}}\|}\times
{\prod _{j=1}^n|L_{i_j}(\ux )|\over |\det(\uL _{i_1}^{tr},\ldots ,\uL _{i_n}
^{tr})|}.\endaligned$$
As with (3), Hadamard's inequality and our bound for $|NS(F)|$ in
Lemma 3 give
$$\prod _{l\neq i_1,\ldots ,i_{n-1}}
{|\det (\uL _l^{tr},\uL _{i_1}^{tr},\ldots ,\uL _{i_{n-1}}^{tr})|
\over\|\uL _l\|\cdot \|\uL _{i_1}\|\cdots \|\uL _{i_{n-1}}\|}\ge |NS(F)|^{1/n!}
\ge
\hofF ^{-b(F)/n!}.$$
Since the discriminant is not zero, each $\uL _i$ occurs in the same
number of\break $n$-tuples in $I(F)$, i.e., $b(\uL _i)=b(F)$ for each
$i$. Hence
$$db(F)=\sum _{i=1}^db(\uL _i)=n|I(F)|.$$
By (2)  then,
$${b(F)\over n!}={n|I(F)|\over dn!}={n\over d}{d\choose n}={d-1\choose n-1}.$$
This proves the case when the discriminant is not zero.
\enddemo

The estimate in Lemma 5 is not so good when $\hofF$ is large in
comparison to $m$ or $\|\ux\|$. In such a situation we will use
the following, which generalizes [S3 Chap.~IV, Lemma 6A].

\nonumproclaim{Lemma 6} Suppose $I(F)$ is not empty 
and $\hofF$ is minimal 
among forms equivalent to $F${\rm .} Suppose further that $F$ does not vanish
at any nonzero integral point{\rm .} Then
for every $\ux\in\br ^n$
there is an $n$\/{\rm -}\/tuple in $I'(F)$
with
$${|F(\ux )|^{n/d}\over \hofF ^{1/d}}
\gg {\prod _{j=1}^n
|L_{i_j} (\ux )|\over |\det (\uL _{i_1}^{tr},\ldots ,\uL _{i_n}^{tr})|}.$$
\endproclaim

\demo{{P}roof} If $F(\ux )=0$ the statement is trivial, so assume otherwise.
By homogeneity of the quantities
$${\prod _{j=1}^n|L_{i_j} (\ux )|\over 
|\det (\uL _{i_1}^{tr},\ldots ,\uL _{i_n}^{tr})|},$$
we may use any 
factorization of $F$. Let $F(\uX )=\prod _{l=1}^kF_l(\uX )=\prod _{i=1}^d
L_i(\uX )$ be the factorization of $F$ in the proof of Lemma 3, and
introduce a new factorization
$F(\uX )=\prod _{i=1}^dL_i'(\uX )$ given by
$$L_i'(\uX )={|F(\ux )|^{1/d}\over |L_i(\ux )|}L_i(\uX )$$ for each $i$. 
By hypothesis,
$${\Cal H}(F\circ T)\ge\hofF =\prod _{i=1}^d\|\uL '_i\|\tag 10$$ 
for any $T\in\gln (\bz )$.

There are $r_1$ real linear factors and
$r_2$ pairs of complex conjugate linear factors of $F$, say.
Arrange the indices so that $\uL '_i\in\br ^n$ for $i\le r_1,$
$\uL '_i\in\bc ^n$ for $r_1<i\le d=r_1+2r_2$ and $\uL '_{i+r_2}=\overline{
\uL '_i}$ for $r_1<i\le r_1+r_2$. Let $\be ^d\subset\br ^{r_1}\oplus 
\bc ^{2r_2}$ 
be the set of
$\ux =(x_1,\ldots ,x_d)$ where $x_{i+r_2}=\overline{x_i}$ for $r_1<i\le r_1+
r_2.$ Then $\be ^d$ is  $d$-dimensional Euclidean space via the usual
hermitian inner product on $\bc ^d$. 

Let $M$ be the $d\times n$ matrix
given by
$$M:=\pmatrix \uL '_1\\ \vdots \\ \uL '_d\endpmatrix =
(\um _1^{tr},\ldots ,\um _n^{tr}).$$
Then $\um _j\in\be ^d$ for all $1\le j\le
n.$ 
Moreover,
$$\|\wedge _{j=1}^n\um _j\| ^2=
\sum _{I'(F)}
 |\det\big ((\uL '_{i_1})^{tr},\ldots ,(\uL '_{i_n})^{tr}\big )|^2,\tag 11$$
where the sum is over all $n$-tuples in $I'(F)$. The interplay between
(10) and (11) which deal with lengths of the rows and columns of $M$, 
respectively,
will be used to get our result.

Let $\lambda _1\le\cdots \le\lambda _n$ be the successive minima of
the $n$-dimensional lattice 
$\Lambda =\oplus _{j=1}^n\bz\um _j\subset\be ^d$ with
respect to the unit ball. Then by Minkowski's theorem,
$$\lambda _1^2\cdots\lambda _n^2 \ll\det (\Lambda )^2=
\|\wedge _{j=1}^n\um _j\|^2.\tag 12$$
We need a lower bound for $\lambda _1^2\cdots \lambda _n^2$.
We first get a lower bound on $\lambda _1^2\cdots\lambda _{n-1}^2$.
We then get a lower bound on $\lambda _n$ and finish
the proof.

Let $\uz _1,\ldots ,\uz _n$ be a basis for $\Lambda$ satisfying $\|\uz _j\|\le
j\lambda _j$ for each $j$. Write
$$MT=(\uz _1^{tr},\ldots ,\uz _n^{tr}),$$
where 
$$T=(\ua _1^{tr},\ldots ,\ua _n^{tr})\in\gln (\bz ),$$ and write
$\uz _j=(z_{j,1},\ldots ,z_{j,d})$ 
for $1\le j\le n$.
Note that $z_{j,i}=L_i'(\ua _j)$ for $j=1,\ldots ,n$ and
$i=1,\ldots ,d$. 
 
Since $F(\ua _1)\neq 0$ by construction, we
have $|F(\ua _1)|\ge 1.$ 
The arithmetic-geometric inequality thus gives
$$\aligned (\lambda _1)^2\ge\|\uz _1\|^{2}&\ge 
d\left (\prod _{j=1}^{d}|z_{1,i_j}|^2\right )^{1/d}\\
&=
d\left (\prod _{j=1}^{d}|L'_{i_j}(\ua _1)|\right )^{2/d}\\
&=
d\left (|F(\ua _1)|\right )^{2/d}\\
&\ge d.\endaligned$$
In particular, 
$$\prod _{j=1}^{n-1}\lambda _j^2\ge\lambda _1^{2(n-1)}\ge 1.\tag 13$$

We need a better bound for $\lambda _n$. For this, we use
another application of the arithmetic-geometric inequality together with (10),
getting
$$\aligned 
n^3\lambda _n^2\ge\sum _{j=1}^n(j\lambda _j)^2\ge
\sum _{j=1}^n\|\uz _j\|^2
&=\sum _{j=1}^n\sum _{i=1}^d|z_{j,i}|^2\\
&=\sum _{i=1}^d\sum _{j=1}^n|z_{j,i}|^2\\
&=\sum _{i=1}^d\|(z_{1,i},\ldots ,z_{n,i})\|^2\\
&\ge d\left (\prod _{i=1}^d\|(z_{1,i},\ldots ,z_{n,i})\|^2\right )^{1/d}\\
&= d\left ({\Cal H}(F\circ T)\right )^{2/d}\\
&\ge d(\hofF )^{2/d}.
\endaligned$$

Our bound for $\lambda _n$ together with the bound
(13) yields
$$\prod _{j=1}^n\lambda _j^2\ge n^{-3}d\hofF^{2/d}.$$
By (11) and (12), we get
$$\sum  _{I'(F)}
|\det\big ((\uL '_{i_1})^{tr},\ldots ,(\uL '_{i_n})^{tr}\big )|^2
\gg \hofF ^{2/d}.$$
There are no more than
${d\choose n}$ summands here by (2).
The largest summand thus satisfies
$$|\det \big ((\uL _{i_1}')^{tr},\ldots ,(\uL _{i_n}')^{tr}\big )|
\gg \hofF ^{1/d}.$$
Finally, we have
$$|\det \big ((\uL _{i_1}')^{tr},\ldots ,(\uL _{i_n}')^{tr}\big )|
={|\det (\uL _{i_1}^{tr},\ldots ,\uL _{i_n}^{tr})|\cdot |F(\ux )|^{n/d}\over
\prod _{j=1}^n
|L_{i_j}(\ux )|}.$$
\enddemo

\vglue-6pt
\section{Auxiliary results}
\vglue-6pt

By Lemmas 5 and 6, any solution to (1) satisfies an inequality of the form
$${\prod _{j=1}^n|L_{i_j}(\ux )|\over|\det (\uL _{i_1}^{tr},\ldots ,
\uL _{i_n}^{tr})|}\ll A,$$
where $A$ is some given bound.
Our goal here is to get information on the solutions to such inequalities.
Specifically, we show that such solutions lie in convex sets. Further,
given bounds for the lengths of such solutions considered, there are
upper bounds for the
number of such convex sets. Lastly, we determine upper bounds for
both the volume and the number of integral points in such convex sets.

\nonumproclaim{Lemma 7} Let $K _1(\uX ),\ldots ,K _n(\uX )\in\bc [\uX ]$ 
be $n$ linearly independent\break linear forms
in $n$ variables{\rm .} Denote the corresponding coefficient vectors
by\break $\uK _1,\ldots ,\uK _n$. Let $A,B,C>0$ with $C>B$ and let $D>1${\rm .} 
Consider the set of $\ux\in\br ^n$ satisfying
$${\prod _{i=1}^n|K_i(\ux )|\over |\det(\uK _1^{tr},\ldots ,\uK _n^{tr}
)|}\le A\tag 14$$
and also $B\le\|\ux\|\le C.$
If $BC^{n-1}\ge D^{n-1}n!n^{n/2}A${\rm ,}
then this set
lies in the union of less than 
$$n^3\left (\log _D\big (BC^{n-1}/n!n^{n/2}A\big )\right )^{n-2} $$
convex sets of the form
$$\gather
\{\uy\in\br ^n\: |K_i'(\uy )|\le a_i\ \text{for $i=1,\ldots ,n$}\}, \\
|\det \big ((\uK _1')^{tr},\ldots ,(\uK _n')^{tr}\big )|=1,\tag 15 \\
\|\uK '_i\| =1\qquad i=1,\ldots ,n,\endgather$$
with 
$$\prod _{i=1}^na_i<D^nn!n^{n/2}{CA\over B}.$$
If $BC^{n-1}<D^{n-1}n!n^{n/2}A,$ then this set lies in the
union of no more than $n!$ convex sets of this form{\rm .}
\endproclaim

\demo{{P}roof} The proof of [S3 Chap.~IV, Lemma 7A] shows that the solutions
to (14) can be partitioned into $n!$ subsets, and for each such subset there
exist pairwise orthogonal linear forms $K '_1(\uX),\ldots ,K ' _n(\uX)$
(these depend on the subset) 
such that all solutions $\ux$ in that subset satisfy
$${\prod _{i=1}^n|K_i'(\ux )|\over |\det\big ((\uK '_1)^{tr},\ldots ,(\uK _n')
^{tr}\big )|}\le n!A.\tag $14'$ $$
After possibly rescaling, we may assume that $\|\uK '_i\|=1$ for each $i$. 
This implies the modulus of the determinant is 1 as well.

Let $\ux$ be a solution to ($14'$) of length at least $B$. 
By Lemma 4, for some $i_0$ (depending on $\ux$, of course) we have
$|K_{i_0}(\ux )|\ge n^{-n/2}B$. This leaves us with
$$ \prod _{i\neq i_0}|K_i'(\ux )|\le {n!n^{n/2}A\over B}.$$

Write
$|K_i'(\ux )|=D^{-n_i}C$ for each $i\neq i_0$. If $\|\ux\|\le C,$ then
$n_i\ge 0$, and by the above estimate $\sum _{i\neq i_0}
n_i\ge\log _D\big (BC^{n-1}/n^{n/2}n!A\big )$.
Let $[\cdot]$ denote the greatest integer
function. Then
$$\sum_{i\neq i_0}[n_i]>\sum _{i\neq i_0}(n_i-1)\ge
\log _D\left ({BC^{n-1}\over D^{n-1}n!n^{n/2}A}\right ).
$$
For the time being, denote the quantity 
$\log _D\left ({BC^{n-1}\over D^{n-1}n!n^{n/2}A}\right )$ by $Q$. If $Q\ge 0$,
then we can find nonnegative integers $z_i\le [n_i]$ for each $i\neq i_0$
that satisfy $\sum _{i\neq i_0}z_i=[Q]$. Further, our solution $\ux$ satisfies
$|K _i'(\ux )|\le D^{-z_i}C$ for $i\neq i_0$ since $-z_i\ge -n_i$.
To make the notation uniform, we set $z_{i_0}=0$, so that
$|K _i'(\ux )|\le D^{-z_i}C$ for all $i$. If $Q<0$, then we simply set
all $z_i=0$.

Summarizing what we have accomplished so far, we see that
the solutions $\ux$ to (14) with $B\le\|\ux\|\le C$
lie in the union of $n!$ subsets, and for each
such subset there are pairwise orthogonal linear forms $K_1'(\uX),
\ldots ,K_n'(\uX )$ with $\|\uK _i'\|=1$ such that all solutions in
that subset lie in convex sets of the form
$$\{\uy\in\br ^n\: |K _i'(\uy )|\le D^{-z_i}C\ \text{for all $i$}\},$$
where the $z_i$'s are nonnegative integers, at least
one of which is 0, satisfying
$$\sum _{i=1}^nz_i=\max \{[Q],0\}.$$
Letting $a_i =D^{-z_i}C$, we have 
$$\prod _{i=1}^na_i\le C^nD^{-[Q]}<C^nD^{1-Q}=D^nn!n^{n/2}{CA\over B}.$$

It remains to estimate the number $n$-tuples $(z_1,\ldots ,z_n)$ which
satisfy the above conditions. Towards that end,
for a nonnegative integer $a$ denote the number of
$n$-tuples $(z_1,\ldots ,z_n)\in\bz ^n$ satisfying $z_i\ge 0$ and
$\sum z_i=a$ by $f(n,a)$. Clearly $f(1,a)=1$. We claim that
$$ f(n,a)\le {(a+n-1)^{n-1}\over (n-1)!}.$$
We see this by induction on $n$. Assuming $n\ge 2$ and our claim is true for
$n-1$,
$$ \aligned \noalign{\vskip5pt} f(n,a)=\sum _{i=0}^af(n-1,a-i)
&\le\sum _{i=0}^a{(a-i+n-2)^{n-2}\over (n-2)!}\\ \noalign{\vskip5pt}
&=\sum _{j=0}^a{(j+n-2)^{n-2}\over (n-2)!}\\ \noalign{\vskip5pt}
&\le {1\over (n-2)!}\int_0^{a+1}(x+n-2)^{n-2}\, dx\\ \noalign{\vskip5pt}
& \le {(a+n-1)^{n-1}\over (n-1)!}.\endaligned$$

Suppose $Q\ge 0$. Then by this claim the number of $n$-tuples 
$(z_1,\ldots ,z_n)$ of
nonnegative integers with $z_{i_0}=0$, say, satisfying $\sum _{i=1}^nz_i=[Q]$
is no greater than $([Q]+n-2)^{n-2}/(n-2)!$. (When $z_{i_0}=0$, we use the
case $n-1$ of our claim.) Taking into account the $n$ 
different possibilities for $i_0$, we see that the total number of
possible $n$-tuples we must consider is no greater than
$n([Q]+n-2)^{n-2}/(n-2)!<n(Q+n-1)^{n-2}/(n-2)!.$ Of course, if $Q<0$ we
have the one $n$-tuple where $z_i=0$ for all $i$.

Now if $Q\ge 0$, we have
$${n(Q+n-1)^{n-2}\over (n-2)!}={n\big (\log _D(BC^{n-1}/n!n^{n/2}A)\big )^{n-2}
\over (n-2)!}.$$
Also, $Q\ge 0$ if and only if $BC^{n-1}\ge D^{n-1}n!n^{n/2}A$. Taking into
account the $n!$ different subsets and using $n\cdot n!/(n-2)!<n^3$ completes
the proof.

\enddemo

We will also use the following variation of Lemma 7, which does away with
the lower bound condition $\|\ux\|\ge B$ at the expense of a higher power
of the logarithmic term in the number of convex sets.

\nonumproclaim{Lemma 7$'$} Let $K _1(\uX ),\ldots ,K _n(\uX )\in\bc [\uX ]$ 
and $\uK _1,\ldots ,\uK _n$ be as in Lemma {\rm 7.} Let $A,C>0$ and $D>1${\rm .}
If $C^{n}\ge D^nn!A${\rm ,}
then the solutions $\ux$ to {\rm (14)} with $\|\ux\|\le C$
lie in the union of less than 
$$n\left (\log _D\big (C^n/n!A\big )\right )^{n-1} $$
convex sets of the form {\rm (15)}
with 
$$\prod _{i=1}^na_i<D^nn!A.$$
If $C^n<D^nn!A,$ then such solutions 
lie in the
union of no more than $n!$ convex sets of this form{\rm .}
\endproclaim

\demo{{P}roof} The proof goes essentially the same way as for Lemma 7.
The difference
is that we do not invoke Lemma 4. Again we have $n!$ subsets where all
solutions in the subset satisfy (14$'$). Let $\ux$ be such a solution
with\break $\|\ux\|\le C$
and write
$|K_i'(\ux )|=D^{-n_i}C$ 
with $n_i\ge 0$ for each $i$. This time we have $\sum _{i=1}^n
n_i\ge\log _D\big (C^n/n!A\big )$, so that
$$\sum_{i=1}^n[n_i]>\log _D\left ({C^n\over D^nn!A}\right ).
$$
This time denote the quantity 
$\log _D\left ({C^n\over D^nn!A}\right )$ by $Q$. As before,
if $Q\ge 0$,
then we can find nonnegative integers $z_i\le [n_i]$ for each $i$
that satisfy $\sum _{i=1}^nz_i=[Q]$ and our solution $\ux$ satisfies
$|K _i'(\ux )|\le D^{-z_i}C$ for all $i$.
If $Q<0$ we  set
all $z_i=0$ again, so that
$$\sum _{i=1}^nz_i=\max \{[Q],0\}.$$
Set $a_i=D^{-z_i}C$ again.

We are now in the same position as with Lemma 7. The difference is that
here we do not  say one of the exponents $z_{i_0}$ is zero, and
$$\prod _{i=1}^na_i\le C^nD^{-[Q]}<C^nD^{1-Q}=D^nn!A.$$
Using the claim in the proof of Lemma 7, the number of $n$-tuples 
$(z_1,\ldots ,z_n)$ of
nonnegative integers  satisfying $\sum _{i=1}^nz_i=[Q]$
is no greater than $$([Q]+n-1)^{n-1}/(n-1)!<(Q+n)^{n-1}/(n-1)!$$ when $Q\ge 0.$ 
Using
$${(Q+n)^{n-1}\over (n-1)!}={\big (\log _D(C^n/n!A)\big )^{n-1}
\over (n-1)!}$$
and $Q\ge 0$ if and only if $C^n\ge D^nn!A$ completes
the proof.
\enddemo

We need estimates for the
number of integer points and also
the volume of the set of all points in convex sets of the form (15).
(When the $K_i(\uX )$s are real linear forms these convex sets are
simply parallelopipeds.)
The following lemmas will provide the needed estimates.

\nonumproclaim{Lemma 8} Let ${\Cal C}\subset\br ^n$ be a convex body {\rm (}\/convex{\rm ,}
closed{\rm ,} bounded and symmetric about the origin\/{\rm )} and let $\Lambda\subset\br ^n$
be a lattice{\rm .} Suppose there are $n$ linearly independent lattice 
points in ${\Cal C}${\rm .} Then there are $\uy _1,\ldots ,\uy _n\in {\Cal C}$ such
that the number of lattice points in ${\Cal C}$ is no greater than
$$3^n2^{n(n-1)/2}{|\det (\uy _1^{tr},\ldots ,\uy _n^{tr})|
\over \det (\Lambda )}.$$
\endproclaim

\demo{{P}roof}
By the homogeneity of the upper bound here we may assume\break
$\Lambda =\bz ^n$. 

The proof is by induction on $n$. If $n=1$, then ${\Cal C}$ is an interval
centered at the origin, say $[-y_1,y_1]$. Since ${\Cal C}$ contains a
nonzero integer point by hypothesis, we have $y_1\ge 1$. Thus, the
number of integer points in ${\Cal C}$ is no greater than $2y_1 +1\le 3y_1.$

Now assume $n>1$ and let $\uz _1,\ldots \uz _n$ be $n$ linearly
independent integer points in ${\Cal C}$. Let $V$ be the span of the
first $n-1$ of them and let $\Lambda ^-=\bz ^n\cap V.$ Then
$\Lambda ^-$ is a primitive sublattice and there is a $\uz _n'\in\bz ^n$
with $\bz ^n=\Lambda ^-\oplus \bz\uz _n'$.

Any integer point $\uz$ in ${\Cal C}$ 
may be written as a sum $\uz =\uz ^-+a\uz _n'$ where $\uz ^-\in\Lambda ^-$
and $a\in\bz$. Further, since $\uz _n$ is an integer point in ${\Cal C}$ but
not in $V$, we see that $a\neq 0$ is possible here. By Cramer's rule
$$|a|={|\det \big ((\uz _1')^{tr},\ldots, (\uz _{n-1}')^{tr},\uz ^{tr}\big )|
\over |\det\big ((\uz _1')^{tr},\ldots, (\uz _n')^{tr}\big )
|}=|\det \big ((\uz _1')^{tr},\ldots ,(\uz _{n-1}')^{tr},\uz ^{tr}\big )|,$$
where $\uz _1',\ldots ,\uz _{n-1}'$ form a basis for $\Lambda ^-$ (so that $\uz _1'\ldots ,\uz _n'$ is a basis for
$\bz ^n$).

For any $a$ we estimate the number of $\uz ^-\in\Lambda ^-$ with
$\uz ^-+a\uz _n'\in {\Cal C}$ as follows. Let $\{\uz _1^-,\ldots ,\uz _N^-\}$
be the set of all such $\uz ^-$. Then the set of differences 
$(\uz _i^- +a\uz _n')-(\uz _1
^- +a\uz _n')$ is a set of $N$ distinct integer points in $\Lambda ^-\cap
2{\Cal C}$ by convexity. Note that $\Lambda ^-$ contains $n-1$
linearly independent lattice points in ${\Cal C}\cap V$, namely $\uz _1,\ldots ,
\uz _{n-1}$. Thus, by the induction hypothesis there are $\uy _1^-,
\ldots ,\uy _{n-1}^-\in 2{\Cal C}\cap V$ such that the number of $\uz ^-\in
\Lambda ^-$ with $\uz ^-+a\uz_n'\in {\Cal C}$ is no greater than
$$3^{n-1}2^{(n-1)(n-2)/2}
{\|\uy _1^-\wedge\cdots\wedge\uy _{n-1}^-\|\over \det (\Lambda ^-)}
.$$
The important thing to note here is the uniformity of this bound; it does not
depend on $a$.

Now let $|a_0|$ be maximal such that there is a $\uz ^-\in\Lambda ^-$ with
$\uz ^-+a_0\uz _n'\in {\Cal C}$ and let $\uy _n$ be this lattice point in $
{\Cal C}$. Let
$\uy _i={1\over 2}\uy _i ^-\in {\Cal C}$ for $i=1,\ldots ,n-1$. We then have
$2|a_0|+1$ possible values of $a$ to consider above, and we now see that
the number of integer points in ${\Cal C}$ is no greater than
$$\aligned & 3^{n-1}2^{(n-1)(n-2)/2}
{\|\uy _1^-\wedge\cdots\wedge\uy _{n-1}^-\|\over \det (\Lambda ^-)}
\big (2|a_0|+1\big )\\
&\qquad \le\ 3^{n-1}2^{(n-1)(n-2)/2} 
{\|\uy _1^-\wedge\cdots\wedge\uy _{n-1}^-\|\over \det (\Lambda ^-)}
3|a_0|\\
&\qquad =\ 3^n2^{n(n-1)/2}{\|\uy _1\wedge\cdots\wedge\uy _{n-1}\|\over \det (\Lambda ^-)}
|\det \big ((\uz _1')^{tr},\ldots ,(\uz _{n-1}')^{tr},(\uy _n)^{tr}\big )|\\
&\qquad =\ 3^n2^{n(n-1)/2}
|\det (\uy _1^{tr},\ldots ,\uy _n^{tr})|.\endaligned$$
\enddemo

Though we do not need it, the proof of Lemma 8 actually shows that
the $\uy _i$'s in ${\Cal C}$ satisfy $2^{n-i}\uy _i\in\Lambda$
as well.

\nonumproclaim{Lemma 9} Let ${\Cal C}$ be a convex body of the form {\rm (15).}
Then either all integral points in
${\Cal C}$ lie in a proper subspace{\rm ,} or the number of such points is
no greater than 
\smallbreak
\centerline{${\displaystyle 3^n2^{n(n-1)/2}n!\prod _{i=1}^na_i.}$}
\smallbreak\noindent 
The volume of ${\Cal C}$ is no greater than
\smallbreak
\centerline{${\displaystyle 2^nn!\prod _{i=1}^na_i.}$}
\endproclaim

\demo{{P}roof} Choose $\uy _1,\ldots ,\uy _n\in {\Cal C}$ with $|\det 
(\uy _1^{tr},
\ldots ,\uy _n^{tr})|$ maximal
(this is clearly possible since ${\Cal C}$ is bounded). Let ${\Cal P}$ be the
region
$${\Cal P}=\{\uy =a_1\uy _1+\cdots +a_n\uy _n\: |a_i|\le 1\ 
\text{for all $i$}\}.$$We claim that ${\Cal P}\supseteq {\Cal C}$. 
Indeed, if there were a $\uy _0\in {\Cal C}\setminus
{\Cal P}$, then without loss of generality $\uy _0=\sum c_i\uy _i$ with
$c_1>1$. But then
$$|\det (\uy _0^{tr},\uy _2^{tr},\ldots ,\uy _n^{tr})|=|\det \big ((c_1\uy _1)
^{tr},\uy _2^{tr},\ldots ,\uy _n^{tr}\big )|
>|\det (\uy _1^{tr},\ldots ,\uy _n^{tr})|,$$
which contradicts the assumption on the $\uy _i$'s. Since the 
volume of ${\Cal P}$ is $2^n|\det (\uy _1^{tr},\ldots ,\uy _n^{tr})|,$ 
we see that
there exist $\uy _1,\ldots ,\uy _n\in{\Cal C}$ with
$$2^n|\det (\uy _1^{tr},\ldots ,\uy _n^{tr})|\ge \text{Vol}({\Cal C}).$$

Finally, if we denote the $n\times n$ matrix with rows $\uK _1',\ldots ,
\uK _n'$ by $T$, then for any 
$\uy _1,\ldots ,\uy _n\in{\Cal C}$ we have
$$\aligned |\det (\uy _1^{tr},\ldots ,\uy _n^{tr})|=
|\det (\uy _1^{tr},\ldots ,\uy _n^{tr})|\times
|\det (T)|&=|\det\big (T\uy _1^{tr},\ldots ,T\uy _n^{tr}\big )|\\
&\le n!\prod _{i=1}^n\max _{1\le j\le n}\{|K_i'(\uy _j )|\}\\
&\le n!\prod _{i=1}^na_i.\endaligned$$
Lemma 9 follows from this estimate, the estimate given above, and
Lemma 8.
\enddemo

\vglue-6pt
\section{The infinite volume case}
\vglue-6pt

This section is devoted entirely  to showing that
the volume $V(F)$ is infinite if $a(F)$
is undefined or at least $d/n$. This is one half of the proposition.
We will also show
that if $a(F)$ is undefined or at least $d/n$, then
(1) has infinitely many integral solutions for $m$ sufficiently large.
Since none of this depends on the particular factorization of $F$ used,
we'll assume that 
$\overline{L _i(\uX )}$ is a factor for all $i$, i.e., the complex
linear factors occur in conjugate pairs.  We break up our argument into
a series of three lemmas.

\nonumproclaim{Lemma 10} If $a(F)$ is undefined or at least $d/n${\rm ,} then there
is a $k<n$ and $k$ coefficient vectors $\uL _{i_1},\ldots ,\uL _{i_k}$
which satisfy the following conditions\/{\rm :}
\smallbreak
\item{\rm 1)}  they are linearly independent\/{\rm ;}
\smallbreak
\item{\rm 2)} there are at least $kd/n$ coefficient vectors $\uL _i$ in their span\/{\rm ;}
\smallbreak
\item{\rm 3)} for all indices $j${\rm ,} if $\overline{\uL _{i_j}}$ is not in the
span of $\uL _{i_1},\ldots ,\uL _{i_j}${\rm ,} then $j<k$ and $\uL _{i_{j+1}}=
\overline {\uL _{i_j}}.$

\endproclaim

\demo{{P}roof} Suppose first that $a(F)$ is undefined. Then by 
Lemma 5 $I(F)$ 
is empty, i.e., the rank of $(\uL _1^{tr},\ldots ,\uL _d^{tr})$ is
less than $n$. Let $k$ be this rank.
Choose an $\uL _{i_1}$. If $k=1$, then all $\uL _i$, in
particular $\overline{\uL _{i_1}}$, are in the span of $\uL _{i_1}$.
If $k>1$, then choose  
an $\uL _{i_2}$ which is linearly independent of $\uL _{i_1}$, with
the stipulation that $\uL _{i_2}=\overline{\uL _{i_1}}$ if this
is a possible choice. Continue on in this fashion, getting
$\uL _{i_1},\ldots ,\uL _{i_k}$. They satisfy conditions 1 and 3 by
construction. There are $d>kd/n$ factors in their span, so condition
2 is satisfied as well.

Now suppose that $a(F)$ is defined and at least $d/n$. Then there is 
an $n$-tuple $(\uL _{i_1},\ldots ,\uL _{i_n})\in J(F)$ and a $j<n$
where $\uL _{i_1},\ldots ,\uL _{i_j}$ have at least $jd/n$ coefficient
vectors in their span
(by the definition of $a(F)$). Let $j_0$ be the least
such index where this is true. By the definition of $J(F)$,
if $\overline{\uL _{i_{j_0}}}$ is in
the span of $\uL _{i_1},\ldots ,\uL _{i_{j_0}}$, then these $j_0$
coefficient vectors satisfy all three conditions above with $k=j_0$.

Suppose 
$\overline{\uL _{i_{j_0}}}$ is not in
the span of $\uL _{i_1},\ldots ,\uL _{i_{j_0}}$. If $j_0=1$, then there
are at least $d/n$ coefficient vectors $\uL _i\not\in\br ^n$ proportional
to $\uL _{i_1}$ and at least $d/n$ additional coefficient vectors 
$\overline{\uL _i}$ proportional to $\overline{\uL _{i_1}}$. In this case
we let $k=2$ and use $\uL _{i_1}$ and $\overline{\uL _{i_1}}$. (Note that
$n>2$ since
$F$ is assumed not to be a power of a positive definite quadratic
form in two variables.) Now suppose $j_0>1$.
Note that condition
3 is still satisfied for all $j<j_0$ by the definition of $J(F).$
Also, by the minimality of $j_0$, there are fewer than $(j_0-1)d/n$
coefficient vectors in the span of $\uL _{i_1},\ldots ,\uL _{i_{j_0-1}}$.
Consider for a moment the collection of $\uL _i$ which are not
in the span of these $j_0-1$ coefficient vectors, but are in the
span of $\uL _{i_1},\ldots ,\uL _{i_{j_0}}.$ We could replace $\uL _{i_{j_0}}$
with any of these and the span would remain the same. 
If $\overline{\uL _i}$ is in
the span of $\uL _{i_1},\ldots ,\uL _{i_{j_0}}$ for one of these
$\uL _i$, then we replace $\uL _{i_{j_0}}$ with $\uL _i$ and let $k=j_0$
as above. If not, then there are more than $(j_0d/n) - (j_0-1)d/n =d/n$ of 
these $\uL _i$, so
there are more than $d/n$ coefficient vectors $\overline{\uL _i}$ which are
not in the span of $\uL _{i_1},\ldots ,\uL _{i_{j_0}}$. This
shows that $j_0d/n$ must be less than $d-(d/n)=(n-1)d/n$, i.e., 
$j_0<n-1$. In this case 
we let $k=j_0 +1<n$ and let $\uL _{i_{k}}=\overline{\uL _{i_{j_0}}}$.
Then conditions 1 and 3 are satisfied. Further, in addition to the
at least $j_0d/n$ coefficient 
vectors in the span of $\uL _{i_1},\ldots ,\uL _{i_{j_0}}$,
we have more than $d/n$ additional coefficient
vectors $\overline{\uL _i}$ in the span of
$\uL _{i_1},\ldots ,\uL _{i_{j_0}},\uL _{i_k}$. This shows that condition
2 holds as well.
\enddemo

\nonumproclaim{Lemma 11} Suppose $a(F)$ is either undefined or at least $d/n${\rm .} Let
$k$ and
$\uL _{i_1},\ldots ,\uL _{i_k}$ be as in Lemma {\rm 10.}
Then
there are linearly independent $\uK _1,\ldots ,\uK _k\in\br ^n$ 
which share the same span as $\uL _{i_1},\ldots ,\uL _{i_k}${\rm .}
\endproclaim

\demo{{P}roof} Suppose $0\le l<k$ and $\uK _1,\ldots ,\uK _l\in\br ^n$ have
been chosen so that their span is equal to the span of $\uL _{i_1},
\ldots ,\uL _{i_l}.$ 

If  
$\overline{\uL _{i_{l+1}}}$ is in the
span of $\uL _{i_1},\ldots ,\uL _{i_{l+1}}$, then it is in the span
of $\uK _1,\ldots , \uK _l,\break \uL _{i_{l+1}}$.  
In this case write 
$$\overline{\uL _{i_{l+1}}}=(a+ib)\uL _{i_{l+1}}+\uz,$$
where $a,b\in\br$ and $\uz\in\bc ^n$ is in the span of $\uK _1,\ldots ,\uK _l.$
Note that both the real and imaginary parts of $\uz$ are in the
span of $\uK _1,\ldots ,\uK _l$ since the $\uK _i$s are real.
A short computation shows that
$$\aligned \Re (\uz )+(a-1)\Re (\uL _{i_{l+1}})&=b\Im (\uL _{i_{l+1}})\\
\Im (\uz ) +b\Re (\uL _{i_{l+1}})&=-(a+1)\Im (\uL _{i_{l+1}}).\endaligned$$
If both $b=0$ and $a=-1,$ then we let $\uK _{l+1}=\Im (\uL _{i_{l+1}}).$
Otherwise we let $\uK _{l+1}=\Re (\uL _{i_{l+1}}).$ In either case
the span of $\uK _1,\ldots ,\uK_l,\uL _{i_{l+1}}$ is equal to
the span of $\uK _1,\ldots ,\uK _l,\uK _{l+1}$.
\smallbreak
If
$\overline{\uL _{i_{l+1}}}$ is not in the
span of $\uL _{i_1},\ldots ,\uL _{i_{l+1}}$, then  $\uL _{i_{l+2}}=
\overline{\uL _{i_{l+1}}}.$ We let
$\uK _{l+1}=\Re (\uL _{i_{l+1}})$ and $\uK _{l+2}=\Im (\uL _{i_{l+1}})$
in this case.
Then the span of 
$\uL _{i_1},\ldots ,\uL _{i_{l+2}}$ is equal to the span of
$\uK _1,\ldots ,\uK_{l+2}$.

Proceeding in this fashion until $l=k$ yields the lemma.

\enddemo

\nonumproclaim{Lemma 12} Suppose $a(F)$ is either undefined or at least
$d/n${\rm .} Let $k$ be as in Lemma {\rm 10.} Then there is an orthonormal basis
$\uK '_1,\ldots ,\uK '_n\in\br ^n$ of $\br ^n$
such that{\rm ,} for all $\ux\in\br ^n$ and $0<a\le b$ satisfying
$$|K'_i(\ux )|\le a\qquad i=1,\ldots ,k$$ and 
$$|K_i'(\ux )|\le b\qquad
i=k+1,\ldots ,n,$$
we have
$$|F(\ux )|^{n/d}\le n^n\hofF ^{n/d}a^kb^{n-k}.$$
Further{\rm ,} $V(F)$ is infinite and $N_F(m)$ is infinite for all $m$
sufficiently large{\rm .}
\endproclaim

\demo{{P}roof} Get $\uK _1,\ldots ,\uK _k$ as in Lemma 11. Let
$\uK _1',\ldots ,\uK _k'$ be an orthonormal basis for their span, and
enlarge this collection to an orthonormal basis $\uK '_1,\ldots ,\uK _n'$
of $\br ^n$. Let $\ux$, $a$ and $b$ be as in the statement of the lemma.
Now at least $kd/n$ of the coefficient vectors $\uL _i$ are
in the span of $\uK '_1,\ldots ,\uK _k'$, and the corresponding factors of
$F(\uX )$  satisfy
$$|L_i(\ux )|\le n\|\uL _i\|\max _{1\le j\le k}\{|K_j'(\ux )|\}=n\|\uL _i\|a.$$
There are no more than $d-kd/n =(n-k)d/n$ factors $L_i(\uX )$ which remain,
and they satisfy 
$$|L_i(\ux )|\le n\|\uL _i\|\max _{1\le j\le n}\{|K_j'(\ux )|\}=n\|\uL _i\|b.$$
Thus,
$$|F(\ux )|=\prod _{i=1}^d|L_i(\ux )|\le n^d\hofF a^{kd/n}b^{(n-k)d/n},$$
and the first part of the lemma is proven.

For $\ux\in\br ^n$ write $\ux =\sum _{i=1}^nx_i\uK '_i$. For any $a\le 1,$
the set of $\ux$ satisfying
$$|x_i|\le\cases  a^{-k/(n-k)}&\text{if $i>k$,}\\
a&\text{if $i\le k$}\endcases$$
is contained in the set of $\ux$ satisfying $|F(\ux )|\le n^d\hofF $
by the first part of the lemma. Letting ${\frak m}$ denote
$\displaystyle{\max _{1\le i\le k}\{|x_i|\}}$ in what follows, we see that
$$\aligned
&\hskip-1in \idotsint\limits_{|x_i|\le 1}\left [\quad\idotsint 
\limits_{|x_j|\le {\frak m}^{-k/(n-k)}}
\prod _{j=k+1}^ndx_j\right ]\prod _{i=1}^kdx_i\\ \noalign{\vskip5pt}
&\qquad = 2^{n-k}
\idotsint\limits_{|x_i|\le 1}
{\frak m}^{-k}
\prod _{i=1}^kdx_i\\ \noalign{\vskip5pt}
&\qquad \ge\idotsint\limits_{\|(x_1,\ldots ,x_k)\|\le 1}
\|(x_1,\ldots ,x_k)\|^{-k}
\prod _{i=1}^kdx_i\\ \noalign{\vskip5pt}
&\qquad =kV(k)\int _0^1r^{-1}dr\\ \noalign{\vskip5pt}
&\qquad =\infty,\endaligned$$
where $V(k)$ denotes the volume of the unit ball in $\br ^k$.
Thus the volume of the set of $\ux\in\br ^n$ with $|F(\ux )|\le n^d\hofF$ is
infinite. By homogeneity, this
shows that $V(F)$ is infinite.

Finally, let $0<a\le b$ satisfy $a^kb^{n-k}=1$. 
Then the parallelopiped
defined by $|K_j'(\ux )|\le a$ for $1\le j\le k$ and $|K_j'(\ux )|\le b$
for $j>k$ has volume $2^n$. By Minkowski's theorem there is a 
nontrivial
integral point in such a parallelopiped. Letting $a\rightarrow 0,$ we
get infinitely many nonzero integral points contained in such 
parallelopipeds. Thus, there are infinitely many integral $\ux$ with
$|F(\ux )|\le n^d \hofF.$
\enddemo

\section{Small solutions}

Let $B_0\ge 1$. 
Any solution 
$\ux \in \br ^n$ to (1) with $\|\ux\|\le B_0$ will be called a {\it small
solution}. We will use 
$B_0=m^{1/(d-a(F))}$ in our proofs of the theorems, but since most of our
estimates up until that point will not require 
``small" to be dependent on $m$, we will leave $B_0$ variable when possible.
In this section we will bound both the volume of all small real solutions
to (1) and the number of small integral solutions, and we will also
also compare the volume of all small solutions with the number of small
integral solutions. As a notational convenience, let
$S_0$ denote the cardinality of the set of small integral solutions
and let $V_0$ denote the volume of all small solutions. 

\nonumproclaim{Lemma 13} Suppose $I(F)$ is not empty{\rm ,}
$\hofF$ is minimal among forms equivalent to
$F${\rm ,} $\hofF >1$ and $F$ has no nontrivial integral
zeros{\rm .} Then
$$V_0\ll m^{n/d}\left (1+{\log B_0\over\log\hofF }\right )^{n-1}$$
and
$$S_0\ll 
m^{n/d}\left (1+{\log B_0\over\log\hofF }\right )^{n-1} +B_0^{n-1}
\left (1+{\log B_0\over\log\hofF }\right )^{n-1}.$$
\endproclaim

\demo{{P}roof} According to Lemma 6, for any solution $\ux\in\br ^n$ to (1)
there is an $n$-tuple in $I'(F)$ with
$$\prod _{j=1}^n{|L_{i_j}(\ux )|\over |\det (\uL _{i_1}^{tr},\ldots ,\uL _{
i_n}^{tr})|}\ll {m^{n/d}\over \hofF ^{1/d}}.$$
Set $A=m^{n/d}/\hofF ^{1/d}$, $C=B_0$ and $D=\hofF^{1/nd}$ in Lemma 7$'$. We
see that the solutions $\ux$ to the above inequality with $\|\ux\|\le C$
lie in 
$$\ll  1+ 
(\log _DC)^{n-1} \ll
\left (1+{\log B_0\over\log\hofF }\right )^{n-1}$$
convex sets of the form (15) with
$$\prod _{j=1}^na_i\ll D^nA= m^{n/d}.$$
By Lemma 9, such a convex set has volume $\ll m^{n/d}.$ There are
no more than ${d\choose n}$ $n$-tuples to consider here by (2), so we
get our bound for $V_0.$

As for $S_0$, we estimate exactly as above. The difference is that our
convex sets may not contain $n$ linearly independent integral points;
they may lie in a proper rational subspace.
So it remains to estimate the number of integral points in these
proper subspaces. By (2) again, there are 
$$\ll \left (1+{\log B_0\over\log\hofF }\right )^{n-1}$$
such subspaces to deal with. We claim that for any proper rational
subspace of $\bq ^n$ of dimension $n'$,
the number of integral points in the subspace
with length at most $B_0$ is $\ll B_0^{n'}.$ Our proof will be complete
once we show this claim.

We prove our claim by induction on $n'$. If $n'=1$ the result is obvious.
Now suppose $W$ is a proper rational subspace of dimension $n'>1.$ 
Let $\Lambda$ be the lattice of integral points in $W$.
If $\Lambda$ doesn't contain $n'$ linearly independent
points of length no more than $B_0,$ then we apply the induction
hypothesis to the proper subspace of $W$ these small lattice points span (and
use $B_0\ge 1$) to show that $\Lambda$ contains $\ll B_0^{n'}$ lattice
points of length at most $B_0$. 

Suppose $\Lambda$ contains $n'$ linearly independent lattice points of length
at most $B_0$. Let $T\in\gln (\br )$ be an orthonormal transformation
taking $W$ to the span of the first $n'$ canonical basis vectors of $\br ^n$.
Let ${\Cal C}\subset \br ^{n'}$ be the set of points of length at most
$B_0$. Since $T(\Lambda )$ is a lattice containing $n'$ linearly independent
lattice points in ${\Cal C}$ and $T$ is orthonormal, Lemma 8 gives
$$|{\Cal C}\cap T(\Lambda )|\ll {B_0^{n'}\over \det (T(\Lambda ))}=
{B_0^{n'}\over \det (\Lambda )}.$$
It is well known that $\det (\Lambda )\ge 1$, so we see that the number
of integral points in $W$ with length at most $B_0$ is $\ll B_0^{n'}$.
Our claim follows by induction, whence our proof of Lemma 13 is complete.
\enddemo

For the purposes of Theorem 3, we need to compare the number of integral
small solutions with the total volume of all small solutions. It proves
convenient here to use the sup norm rather than the Euclidean norm.
So let $V_0'$ denote the volume of all solutions to (1) with sup norm
at most $B_0$, and similarly for $S_0'$.

\nonumproclaim{Lemma 14} With the notation above{\rm ,} we have
$$|S_0'-V_0'|\le dn(2B_0+1)^{n-1}.$$
\endproclaim

\demo{{P}roof} Let $\mu$ denote the usual Lebesgue measure on $\br$ and let
$\nu$ denote the $\sigma$-finite measure gotten from the characteristic
function of $\bz$, that is, $\nu (E)$ is the number of integer points in
the set $E$ for any Borel set $E\subseteq \br$. Let $\chi$
be the characteristic function of the set
$$\{\uy\in\br ^n\: |F(\uy )|\le m\ \text{and $|y_i|\le B_0$ for all $i$}\}.$$
What we want to do here is estimate the difference between the integrals
of $\chi$ with respect to the product measures $\mu ^n$ and $\nu ^n$. The
lemma follows from the case $I=\{1,\ldots ,n\}$ of the following claim:

For any nonempty subset $I\subseteq \{1,\ldots ,n\}$ 
and fixed values $y_i\in\br$ for
$i\not\in I,$ we have
$$\multline
\left |\int\cdots\int\chi (y_1,\ldots ,y_n)\prod _{i\in I}d\mu (y_i)-
\int\cdots\int\chi (y_1,\ldots ,y_n)\prod _{i\in I}d\nu (y_i)\right |\\
\le
d|I|(2B_0+1)^{|I|-1},\endmultline $$
where $|I|$ denotes the cardinality of $I$.
The major point of this estimate is that it is independent of the
particular choices of $y_i\in\br$ for $i\not\in I$. We prove this claim
(and whence Lemma 14) by induction on the cardinality of $I$. 

Suppose that $I=\{i_0\}$ and $y_i\in\br$ are fixed for $i\neq i_0$. Then
$$F(y_1,\ldots ,Y_{i_0},  \ldots ,y_n)\in\br [Y_{i_0}]$$ is a polynomial in
one variable  of degree no greater than $d$. This implies that the set
$$E= \{y_{i_0}\in\br \: |F(y_1,\ldots ,y_n)|\le m\ \text{and $|y_i|\le B_0$
for all $i$}\}$$
is a (possibly empty) union of no more than $d$ nonintersecting closed
intervals. Now 
$$\int\chi (y_1,\ldots ,y_n)d\mu (y_{i_0})=\int _Ed\mu (y_{i_0})$$
and
similarly for the $\nu$ measure. Further, the difference between the
length of a closed interval and the number of integer values therein
is between $-1$ and~$1$. This shows the case $|I|=1$ of the claim.

Now suppose $|I|>1$. We will use the induction hypothesis twice and
the Fubini-Tonelli theorem to show the claim holds for $I$. Choose $i_0\in I$.
Then by the Fubini-Tonelli theorem and the triangle inequality

{\ninepoint $$\multline 
\left |\idotsint\chi (y_1,\ldots ,y_n)\prod _{i\in I}d\mu (y_i)-
\idotsint\chi (y_1,\ldots ,y_n) \prod_{i\in I}d\nu (y_i)\right |\\
\le
\left |\int
\left [\int\cdots \int \chi (y_1,\ldots ,y_n)
\!\prod _{i\in I, i\neq i_0}\! d\mu (y_i)-
\int\cdots\int\chi (y_1,\ldots ,y_n)
\!\prod _{i\in I, i\neq i_0}\!d\nu (y_i)\right] \, d\mu (y_{i_0})\right |\\
+
\left |\int\cdots\int
\left [\int\chi (y_1,\ldots ,y_n)\, d\mu (y_{i_0})
-\int\chi (y_1,\ldots ,y_n)\, d\nu (y_{i_0})\right] \prod _{i\in I, i\neq i_0}d\nu (y_i)\right |.
\endmultline$$}
\eject
\noindent 
Using the induction hypothesis on $I\setminus\{i_0\}$ and the fact
that $\chi$ is the characteristic function of a set contained in the
cube $\{\uy\in\br ^n\: |y_i|\le B_0\}$ gives

{\ninepoint $$\multline
\left |\int\left [\int\cdots\int \chi (y_1,\ldots ,y_n)
\!\prod _{i\in I, i\neq i_0}\! d\mu (y_i)-
\int\cdots\int\chi (y_1,\ldots ,y_n)
\!\prod _{i\in I, i\neq i_0}\! d\nu (y_i)\right ]\, d\mu (y_{i_0})\right |\\
\le
\int \left |\int\cdots\int\chi (y_1,\ldots ,y_n)
\prod _{i\in I, i\neq i_0}d\mu (y_i)-
\int\cdots\int\chi (y_1,\ldots ,y_n)
\prod _{i\in I, i\neq i_0}d\nu (y_i)\right |\, d\mu (y_{i_0})\\
\le\int _{[-B_0,B_0]}d(|I|-1)(2B_0+1)^{|I|-2}\, d\mu (y_{i_0})\\
=2B_0d(|I|-1)(2B_0+1)^{|I|-2}< d(|I|-1)(2B_0+1)^{|I|-1}.\endmultline$$}

\noindent 
Similarly,

{\ninepoint $$\multline
\left |\int\cdots\int
\left [\int\chi (y_1,\ldots ,y_n)\, d\mu (y_{i_0})
-\int\chi (y_1,\ldots ,y_n)\, d\nu (y_{i_0})\right] \prod _{i\in I, i\neq i_0}d\nu (y_i)\right |\\
\le
\int\cdots\int\left |\int\chi (y_1,\ldots ,y_n)\, d\mu (y_{i_0})-
\int\chi (y_1,\ldots ,y_n)\, d\nu (y_{i_0})\right |
\prod _{i\in I, i\neq i_0}d\nu (y_i)\\
\le\int _{[-B_0,B_0]}\cdots\int _{[-B_0,B_0]}d
\prod _{i\in I, i\neq i_0}d\nu (y_i)\\
=d(2[B_0]+1)^{|I|-1}\le d(2B_0+1)^{|I|-1}.\endmultline$$}

\noindent 
Adding these two estimates together finishes our proof of the claim.
\enddemo

\section{Estimating large solutions}

Throughout this section we will assume that $a(F)$ is defined and
less than $d/n$ (this forces $d>n$). It is appropriate at this time to note some
inequalities involving $a(F)$ and $c(F)$ under this assumption.
By definition, $ka(F)\in\bz$ for some $k<n$, so that
$kna(F)\le kd -1$ and 
$$1\le a(F)\le {d\over n}-{1\over n(n-1)}.\tag 16$$ 
Using this, we get 
$$n-d\le {na(F)-d\over a(F)}\le {1\over 1-n}.\tag 17$$
If the discriminant of $F$ is not zero, then 
$$1\le {d-1\choose n-1}-1=c(F).$$
If the discriminant of $F$ is zero, then by (2), (16) and (17)
$$\aligned c(F)&= {b(F)\over n!}\times {d-(n-1)a(F)\over a(F)}-{1\over a(F)
}\\
&< |I'(F)|(d-n+1)\\
&\le {d\choose n}(d-n+1).\endaligned$$
Here we also used $b(F)/n!\le |I'(F)|,$ which is clear from the definitions.
Using $b(F)/n!\ge 1$ (which is also clear from the definitions), and $a(F)<d/n$
gives
$$\aligned c(F)&= {b(F)\over n!}\times {d-(n-1)a(F)\over a(F)}-{1\over a(F)
}\\
&\ge{d-na(F)\over a(F)}+{a(F)-1\over a(F)}\\
&>{(d-n)\over d}.\endaligned$$
Thus,
$${(d-n)\over d}\le c(F)\le {d\choose n}(d-n+1).\tag 18$$

For indices $l\ge 0$ let $B_l=e^lB_0$ and $C_l=e^{l+1}B_0$.  Let
$$A_0=m^{1/a(F)}B_0^{(na(F)-d)/a(F)}\hofF ^{c(F)}$$ and  for $l\ge 0$ let
$A_l=e^{(na(F)-d)l/a(F)}A_0$. Recall that $m,B_0\ge 1$ by hypothesis and
$\hofF\ge 1$ by Lemma 3. By (16), (17), 
and (18)
$$A_l=e^{(na(F)-d)l/a(F)}A_0
\ge B_0^{n-d}e^{l(n-d)}.\tag 19$$
Let $V_{l+1}$
denote the total volume of the set of solutions $\ux\in\br ^n$ to (1)
with $B_l\le\|\ux\|\le C_l$.

\nonumproclaim{Lemma 15} If $I(F)$ is not empty and  $a(F)<d/n${\rm ,} then
$$\sum _{l=1}^{\infty}V_l\ll
\hofF ^{c(F)}m^{1/a(F)}B_0^{(na(F)-d)/a(F)}(1+\log B_0)^{n-2}.$$
\endproclaim

\demo{{P}roof of the proposition} Set $m=B_0=1$. Clearly $V_0\ll 1$. By
Lemma 15
$\sum _{l=1}^{\infty}V_l<\infty$ whenever $a(F)$ is
defined and less than $d/n$. This together with Lemma 12 
proves the proposition.
\enddemo

\demo{{P}roof of Lemma {\rm 15}} 
By Lemma 5, for any solution $\ux\in\br ^n$ to (1) with $B_l\le\|\ux\|$
there is an $n$-tuple in $I'(F)$ with
$$\align {\prod _{j=1}^n|L_{i_j}(\ux )|\over|\det (\uL _{i_1}^{tr},\ldots ,
\uL _{i_n}^{tr})|}&\ll \left ({|F(\ux )|\over
\|\ux \|^{d-na(F)}}\right )^{1/a(F)}\hofF^{c(F)}\tag 20\\
&\le\left ({m\over B_{l}^{d-na(F)}}\right )^{1/a(F)}\hofF ^{c(F)}\\
&=m^{1/a(F)}B_0^{(na(F)-d)/a(F)}e^{(na(F)-d)l/a(F)}\hofF ^{c(F)}\\
&=A_l.\endalign$$ 

We will estimate using  Lemma 7. We have
$$\multline   \max \left \{n!,n^3\big (\log (B_lC_l^{n-1}/n!n^{n/2}A_l)\big )^
{n-2}\right \}\\
 \le\max \left \{n!,n^3\big (\log (B_0^de^{n(l+1)}/e^{l(n-d)}
)\big )^{n-2}\right \}\\
 \ll (1+l+\log B_0)^{n-2}
 \endmultline \tag 21$$
by (19). Setting $A=A_l,\ B=B_l,\ C=C_l$ and $D=e$ in Lemma 7,
we see by (21) that the solutions $\ux$ to (20) with $B_l\le\|\ux\|\le C_l$
are contained in $\ll (1+l+\log B_0)^{n-2}$ convex sets of
the form (15) with
$$\prod _{i=1}^na_i\ll {C_lA_l\over B_l}\ll A_l\le e^{l/(1-n)}
A_0$$
by (17).
According to Lemma 9, the volume of such a convex set is\break 
$\ll e^{l/(1-n)}A_0.$
Taking into account the total number of possible $n$-tuples in $I'(F)$
using (2), we find that
$$V_{l+1}\ll e^{l/(1-n)}A_0(1+l+\log B_0)^{n-2}\le {(1+l)^{n-2}\over
(e^{1/(n-1)})^l}A_0
(1+\log B_0)^{n-2}.$$

We thus have
$$\aligned\noalign{\vskip5pt} \sum _{l=0}^{\infty}V_{l+1}
&\ll A_0(1+\log B_0)^{n-2}
\sum _{l=0}^{\infty}{(l+1)^{n-2}\over (e^{1/(n-1)})^l}\\ \noalign{\vskip5pt}
&\ll \hofF ^{c(F)}m^{1/a(F)}B_0^{(na(F)-d)/a(F)}(1+\log B_0)^{n-2}.\\
\noalign{\vskip5pt}
\endaligned$$
\enddemo

When estimating the number of integer solutions to (1) of length
greater than $B_0$, we proceed very much as in the proof of Lemma 15.
However, since we are counting integer solutions as opposed to
estimating volumes, we must also account for the possibility
that all solutions in a given convex set of the form (15)
lie in a proper subspace, so that Lemma 9 cannot be used in a manner
similar to our use of it in the proof above. 
Our goal is to reach the point where we may
estimate the remaining (extremely large) integer solutions using
a quantitative version of the subspace theorem.

\nonumproclaim{Lemma 16} Suppose $I(F)$ is not empty and $a(F)<d/n.$ 
Then 
the integral solutions $\ux$ to {\rm (1)} with $B_0\le \|\ux\|$
lie in the union of a set
of cardinality $S$ satisfying
$$S
\ll m^{1/a(F)}B_0^{(na(F)-d)/a(F)}(1+\log B_0)^{n-2}\hofF ^{c(F)}$$
and 
$$\ll \big (1+\log m+\log\hofF \big )\big (1+\log m+\log\hofF+\log B_0 \big )
^{n-2}$$
proper rational subspaces{\rm .}
\endproclaim

\demo{{P}roof} Exactly as in the proof of Lemma 15, any integral solution
$\ux$ to (1) with $B_l\le \|\ux\|\le C_l$ satisfies (20) for some
$n$-tuple in $I'(F)$. We apply
Lemma 7 again, getting the same convex sets of the form (15) as in the
proof of Lemma 15. When
those sets contain $n$ linearly independent lattice points, we
estimate the number of such points using Lemma 9 exactly as we 
estimated the $V_{l+1}$
above. These points make up the set of cardinality $S$.

By (21), our solutions $\ux$ to (20)
with $B_l\le\|\ux\|\le C_l$ lie in the union of 
$\ll (l+1+\log B_0)^{n-2}$
convex sets of the form (15). Taking into account the different
possible $n$-tuples, we see that those solutions $\ux$ with
$B_0\le\|\ux\|\le C_{l}$ not already accounted for in $S$ lie
in the union of
$\ll (l+1)(l+1+\log B_0)^{n-2}$
proper rational subspaces. 
We need to determine
how large $l$ should be so that solutions of length at least $C_l$
can be dealt with using the subspace theorem.

If 
$${l+1\over 2(n-1)}\ge \log m +{d\choose n}(d-n+1)\log (\hofF),$$
then by Lemma 2, (16), (17) and (18) we have
$$\aligned C_l^{(d-na(F))/2a(F)}\ge C_l^{1/2(n-1)}&\ge e^{(l+1)/2(n-1)}\\
&\ge m\hofF ^{{d\choose n}(d-n+1)}\\
&\ge m^{1/a(F)}\hofF ^{c(F)}.\endaligned$$
By Lemma 5 and (17), for any solution
$\ux$ to (1) with $\|\ux\|\ge C_l$ there is an $n$-tuple
in $I'(F)$ satisfying
$$\aligned{\prod _{j=1}^n|L_{i_j}(\ux )|\over|\det (\uL _{i_1}^{tr},\ldots ,
\uL _{i_n}^{tr})|}&\ll 
\left ({|F(\ux )|\over \|\ux\| ^{d-na(F)}}\right )
^{1/a(F)}\hofF ^{c(F)}\\
&\le {1\over \|\ux\| ^{(d-na(F))/2a(F)}}
\left ({m\over C_l^{(d-na(F))/2}}\right )^{1/a(F)}\hofF ^{c(F)}\\
&\le \|\ux\| ^{-1/2(n-1)}
\left ({m\over C_l^{(d-na(F))/2}}\right )^{1/a(F)}\hofF ^{c(F)}.\endaligned$$
Let $l_0$ be least such that $\|\ux\|\ge C_{l_0}$ implies that
$${\prod _{j=1}^n|L_{i_j}(\ux )|\over|\det (\uL _{i_1}^{tr},\ldots ,
\uL _{i_n}^{tr})|}<
\|\ux\| ^{-1/2(n-1)}.$$
By what we showed above, 
$l_0\ll 1+\log m+\log\hofF .$ Let $l_1$ be the least such that
$C_{l_1}\ge m^{1/d}C_{l_0}$ and $C_{l_1}\ge m^{1/d}\hofF$. 
Then $l_1\ll 1+\log m+\log\hofF$, too.

The integral solutions $\ux$ to (1) with $B_0\le\|\ux\|\le C_{l_1}$
either lie in our set of cardinality $S$ or
$$\multline \ll l_1(l_1+1+\log B_0)^{n-2}\\
\ll\big (1+\log m+\log\hofF \big )
\big (1+\log m+\log \hofF +\log B_0)^{n-2}\endmultline $$
proper rational subspaces. 
Since the solutions to $F(\ux )=0$ lie in no more than $d$ proper 
subspaces, we
restrict ourselves for what remains to integral solutions $\ux$
to (1) with $|F(\ux )|\ge 1$ and $\|\ux\|\ge C_{l_1}$.
Let $\ux$ be such a solution and
write $\ux =g\ux '$ for some primitive integer point $\ux '$ and some
integer $g\ge 1$. 
By the
homogeneity of $F$, 
$$m\ge |F(\ux )|=|F(g\ux ')|=g^d|F(\ux ')|\ge g^d,$$
so that $g\le m^{1/d}$. 
Thus, 
$$\|\ux '\|\ge m^{-1/d}\|\ux\|\ge m^{-1/d}C_{l_1}\ge \max \{C_{l_0},\hofF\}$$
and $\ux '$ is a primitive solution to (1). 

By the definition of $l_0$, we have
$${\prod _{j=1}^n|L_{i_j}(\ux ')|\over|\det (\uL _{i_1}^{tr},\ldots ,
\uL _{i_n}^{tr})|}<\|\ux '\| ^{-1/2(n-1)}$$
for some $n$-tuple in $I'(F)$.
By Lemma 2 we may assume each $\uL _{i_j}$ here is defined over
a number field of degree at most $d$ and has field height at most
$\hofF\le\|\ux '\|.$ 
By a version of the
quantitative subspace theorem due to Evertse [E2, Corollary], the set of
such primitive integral $\ux '$ 
lies in the union of $\ll 1$ proper subspaces. Taking into account the
number of possible $n$-tuples using (2), we see that the integral
solutions $\ux$ to (1) with $\|\ux\|\ge C_{l_1}$ lie in $\ll 1$ proper
rational subspaces. This completes the proof.
\enddemo

\section{Proof of the theorems}

As remarked above, to prove our theorems we set
$B_0=m^{1/(d-a(F))}$, giving
$$\align m^{1/a(F)}B_0^{(na(F)-d)/a(F)}&=m^{1/a(F)}m^{(na(F)-d)/a(F)(d-a(F))}
\tag 22\\
&=m^{(d-a(F)+na(F)-d)/a(F)(d-a(F))}\\
&=m^{(n-1)/(d-a(F))}.\endalign  $$

\demo{{P}roof of Theorem $1$} By the proposition{\rm ,} it suffices to prove that
$V(F)\ll 1$ when $I(F)$ is not empty and $a(F)<d/n${\rm .} Moreover{\rm ,} by
homogeneity we need only show that $m^{n/d}V(F)\ll m^{n/d}$ for
some positive $m${\rm .} We may assume $\hofF$ is minimal among forms
equivalent to $F$ since $V(F)$ is invariant under equivalence. 

Suppose first that $\hofF =1$. In this case we set $m=1$, too. Clearly
$V_0$ is no larger than the volume of the unit ball in $\br ^n$.
By (22) and Lemma 15 we have
$$\sum _{l=1}^{\infty}V_l\ll 1;$$
thus,
$$V(F)=\sum _{l=0}^{\infty}V_l\ll 1.$$

Now suppose $\hofF >1.$
By (16) we have 
$${n-1\over d-a(F)}\le {n-1\over d-{d\over n}+{1\over n(n-1)}}={n\over 
d+{1\over (n-1)^2}}.\tag 23$$
Choose $m$ so that
$$\hofF ^{c(F)}m^{(n-1)/(d-a(F))}= m^{{n\over d+1/2(n-1)^2}}.$$
Then $\log m\gg\ll\log\hofF$ and 
$$\hofF ^{c(F)}m^{(n-1)/(d-a(F))}(1+\log m )^{n-2}\ll m^{n/d}$$
by (18) and (23).
By Lemma 13,
$$\aligned V_0&\ll m^{n/d}\left (1+{\log B_0\over \log\hofF}\right )^{n-1}\\
&\ll m^{n/d}\left (1+{\log m\over\log\hofF }\right )^{n-1}\\
&\ll m^{n/d}.\endaligned$$
By Lemma 15 and (22),
$$\multline \sum _{l=1}^{\infty} V_l\ll \hofF ^{c(F)}
m^{1/a(F)}B_0^{(na(F)-d)/a(F)}\\ =
\hofF ^{c(F)}m^{(n-1)/(d-a(F))}(1+\log m)^{n-2}\ll m^{n/d}.\endmultline $$
Thus,
$$m^{n/d}V(F)=\sum _{l=0}^{\infty}V_l\ll m^{n/d}.$$
\enddemo

\demo{{P}roof of Theorem {\rm 2}}
Suppose $W$ is a proper rational subspace of $\br ^n$ of dimension $n'$.
Then there is
a $T\in\gln (\bz )$  with
$$T:W\cap\bz ^n\rightarrow \{(z_1,\ldots ,z_n)\in\bz ^n\:
z_i=0\ \text{for $i>n'$}\}.$$
Then $G:=F\circ T^{-1}$ is an equivalent form, and $F$ restricted to $W$ is
equivalent to $G$ restricted to $\br ^{n'}$.
In this manner, we see that considering integral 
solutions to (1) for $F$ restricted
to a proper rational subspace is equivalent to considering integral solutions to
(1) for a form in fewer variables. With this in mind, we will prove that
$N_F(m)\ll m^{n/d}$ when $F$ is of finite type by induction on $n$. But
we first deal with the simpler case when $F$ is not of finite type. 

Suppose that $F$ is not of finite type. Then there is some
nontrivial subspace $W$ defined over $\bq$
where the volume of solutions to
(1) in $W$ is infinite. Let $n'\ge 1$ be the dimension of $W$. If $n'=1,$
then $F$ vanishes on $W$.
Trivially $N_F(m)$ is infinite for all $m$ in this case. Suppose
$n'>1$ and get a
form $F'(\uX )\in\bz [\uX ]$  in $n'$ variables where the $\ux$ in $W$ are
in one-to-one correspondence with $\ux '\in\br ^{n'}$ via a $T\in\gln (\bz )$
with $F(\ux ) =F'(\ux ')$ as above. Since $V(F')$ is infinite
by hypothesis, the proposition shows that $a(F')$ is either undefined or
at least $d/n'$. Lemma 12 shows that $N_{F'}(m)$ is infinite for
all sufficiently large $m$. Thus, there are infinitely many solutions
$\ux\in W\cap
\bz ^n$ to (1) for all sufficiently large $m$.
This shows that $N_F(m)$ is infinite for all sufficiently large $m$ when
$F$ is not of finite type.

Now suppose $F$ is of finite type. Interestingly, our argument for the first
step in the induction where $n=2$ is the same as our argument for $n>2$
using the induction hypothesis. Rather than present the same argument
twice, then, we will simply assume that $n\ge 2$ and that the number
of integral solutions to (1) restricted to a proper subspace of
dimension $n'<n$ is $\ll m^{n'/d}$. 
The number of solutions to (1) restricted
to any proper 1-dimensional rational subspace is $\ll m^{1/d}$, since $F$
is not identically $0$ on such a subspace, so our assumption in the
case $n=2$ is correct. Finally, 
without loss of generality
we may assume $\hofF$ is minimal among forms equivalent to $F$.

By the proposition, $I(F)$ is not empty
and $a(F)<d/n$. 
Suppose first that
$$\hofF ^{c(F)}m^{(n-1)/(d-a(F))}\le m^{{n\over d+1/2(n-1)^2}}.$$
Then (18) and (23) show that $\log\hofF\ll\log m$, and we also have
$$\hofF ^{c(F)}m^{(n-1)/(d-a(F))}(1+\log m)^{n-2}\ll m^{n/d}.$$
By Theorem 1, $V_0'\le m^{n/d}V(F)\ll m^{n/d}$, so (22), (23) and Lemma 14
give $S_0'\ll m^{n/d}.$ Further, (22), (23), and Lemma 16 
show that the integral solutions
of length at least $B_0$ lie in the union of a set of cardinality
$\ll m^{n/d}$ and $\ll (1+\log m)^{n-1}$ proper subspaces. 
By the induction hypothesis (or the trivial 1-dimensional case when $n=2$),
these proper subspaces
contribute 
$$\ll m^{(n-1)/d}(1+\log m)^{n-1}\ll m^{n/d}$$ integral solutions. 
So
$N_F(m)\ll m^{n/d}.$
 
Now suppose
$$\hofF ^{c(F)}m^{(n-1)/(d-a(F))}\ge m^{{n\over d+1/2(n-1)^2}} ,$$
so that $\log\hofF\gg 1+\log m$ by (18) and (23).
Let $C_{l_1}$ be as in the proof of Lemma 16. As shown in the proof
of Lemma 16, 
$l_1\ll (1+\log m+\log \hofF)$.

By Lemma 6, if $\ux$ is a solution
to (1), then there is a $n$-tuple in $I'(F)$ such that
$${\prod _{j=1}^n|L_{i_j}(\ux )|\over|\det (\uL _{i_1}^{tr},
\ldots ,\uL _{i_n}^{tr})|}\ll m^{n/d}\hofF ^{-1/d}.$$
We use Lemma $7'$ with $A=m^{n/d}\hofF ^{-1/d}$, $C=C_{l_1}$ 
and $D=\hofF ^{1/nd}$. 
We have
$$\log _D\left ({C^n\over n!A}\right )\ll {\log C_{l_1}\over\log\hofF }\ll
{\log m+l_1\over\log\hofF }\ll 1,$$
so the set of all such $\ux$ with $\|\ux\|\le C_{l_1}$
lie in $\ll 1$ convex sets of the form (15) with
$$\prod _{i=1}^na_i\ll D^nA=m^{n/d}.$$
By Lemma 9, if such a set contains $n$ linearly independent integral points,
it contains $\ll m^{n/d}$ of them. Taking into account the number of
possible $n$-tuples via (2), we see that the integral solutions $\ux$
to (1) with $\|\ux \|\le C_{l_1}$ lie in the union of $\ll 1$ proper 
rational subspaces and a set of cardinality $\ll m^{n/d}$. 
As shown in the proof of Lemma 16, all integral solutions $\ux$ to (1)
with $\|\ux\|\ge C_{l_1}$ lie in $\ll 1$ proper subspaces.
By the induction hypothesis (or the trivial 1-dimensional case if $n=2$),
all our proper subspaces contain
$\ll m^{(n-1)/d}$ integral solutions total.
So $N_F(m)\ll m^{n/d}$.
\enddemo

\phantom{SW}

\demo{{P}roof of Theorem {\rm 3}} By Lemma 15 and (21)
we have
$$m^{n/d}V(F)-V_0'\le \sum _{l=1}^{\infty}V_l\ll \hofF ^{c(F)}m^{(n-1)/
(d-a(F))}(1+\log m)^{n-2}.$$
By Lemma 14, we get
$$|S_0'-m^{n/d}V(F)|
\ll \hofF ^{c(F)}m^{(n-1)/
(d-a(F))}(1+\log m)^{n-2}.$$
As we saw in the proof of Theorem 2,
the number
of integral solutions to (1) restricted to any proper subspace is $\ll
m^{(n-1)/d}.$ By Lemma 16 and (21) then, the number of integral solutions
to (1) with length at least $B_0$ is 
$$\ll \hofF ^{c(F)}m^{(n-1)/(d-a(F))}(1+\log m)^{n-2}.$$
From this, we get
$$N_F(m)-S_0'\ll
\hofF ^{c(F)}m^{(n-1)/(d-a(F))}(1+\log m)^{n-2}.$$
Theorem 3 follows.
\enddemo

\AuthorRefNames [BT]
\references 

[B] \name{M. Bean}, An isoperimetric inequality for the area of
plane regions defined by binary forms, {\it Compositio Math.}\/ {\bf 92} (1994), 115--131.

[BT] \name{M. Bean} and \name{J.\ L. Thunder}, Isoperimetric inequalities
for volumes associated with decomposable forms, {\it J. London Math.\ Soc.}\/ {\bf 54} (1996), 39--49.

[E1] \name{J.-H. Evertse}, On equations in S-units and the
Thue-Mahler equation, {\it Invent.\ Math.} {\bf 75} (1984), 561--584.

[E2] \bibline,  An improvement of the
quantitative subspace theorem, {\it Compositio Math.}\/ {\bf 101} (1996), 225--311.

[E3] \bibline,  On the norm form inequality $|F(\ux )|\le h$, {\it Publ.\ Math.\ Debrecen} {\bf 56} (2000), 337--374.

[M] \name{K. Mahler}, Zur Approximation algebraischer Zahlen III, {\it Acta Math.}\/ {\bf 62} (1934), 91--166.

[R] \name{K. Ramachandra}, A lattice-point problem for norm forms
in several variables, {\it J. Number Theory} {\bf 1} (1969), 534--555.

[S1] \name{W.\ M. Schmidt}, Norm form equations, {\it Ann.\ of Math.}\/ {\bf 96} (1972), 526-551.

[S2] \bibline, The number of solutions of norm form equations, {\it Trans.\ Amer.\ Math.\ Soc.}\/ {\bf 317} (1989),
197--227.

[S3] \name{W.\ M. Schmidt}, {\it Diophantine Approximations and Diophantine
Equations}, {\it Lecture}\break {\it Notes in Mathematics} {\bf 1467},  Springer-Verlag, New York, 1991.

[T] \name{A. Thue}, \"Uber Ann\"aherungwerte algebraischer Zahlen, {\it J. Reine Angew.\ Math.}\/ {\bf 135} (1909),
284--305.

\endreferences
\bye